\documentclass[12pt,fleqn, a4]{article}

\usepackage[french]{babel}
\usepackage{inputenc}
\usepackage[T1]{fontenc}

\usepackage{amssymb}
\usepackage{latexsym}
\usepackage{graphics}

\usepackage{graphics}
\usepackage{color}
\newcommand{\colval}{0.3}
\definecolor{colone}{gray}{\colval}

\newcommand{\dcb}{\begin{array}{lll}}
\newcommand{\dce}{\end{array}}
\newcommand{\ebe}{\begin{enumerate}\setlength{\baselineskip}{13pt}\setlength{\parskip}{5pt}}
\newcommand{\dbe}{\end{enumerate}}

\newcommand{\ibegin}{\begin{itemize}\setlength{\baselineskip}{19pt}\setlength{\parskip}{7pt}}
\newcommand{\iend}{\end{itemize}}
\newcommand{\ok}{\rule{4pt}{6pt}}
\newcommand{\desb}{\begin{description}}
\newcommand{\dese}{\end{description}}

\newtheorem{Thm}{Theorem}[section]
\newtheorem {Cor}[Thm]{Corollary}
\newtheorem {definition}{Definition}[section]
\newtheorem {pro}{Proposition}[Thm]
\newtheorem {Lemma}[Thm]{Lemma}
\newtheorem {rem}{Remark}[section]

\newtheorem {assumption}{Assumption}[section]

\newcommand {\bd}{\begin{definition}}
\newcommand {\ed}{\end{definition}}
\newcommand {\bpro}{\begin{pro}}
\newcommand {\epro}{\end{pro}}
\newcommand {\bl}{\begin{Lemma}}
\newcommand {\el}{\end{Lemma}}
\newcommand {\bcor}{\begin{Cor}}
\newcommand {\ecor}{\end{Cor}}
\newcommand {\brem }{\begin{rem} \rm }
\newcommand {\erem }{\end{rem}}
\newcommand{\bethe}{\begin{Thm}}
\newcommand{\ethe}{\end{Thm}}
\newcommand {\bassumption}{\begin{assumption}}
\newcommand {\eassumption}{\end{assumption}}

\newcommand{\mrt}{{$\mathfrak{M}$rp}}
\newcommand{\transp}{{^\top\!}}

\def \ind{1\!\!1}

\begin{document}

\begin{center}
\Large
Drift operator in a market affected by the expansion of information flow : a case study
\end{center}

\begin{center}
Shiqi Song

{\footnotesize Laboratoire Analyse et Probabilités\\
Université d'Evry Val D'Essonne, France\\
shiqi.song@univ-evry.fr}
\end{center}

\

\

\begin{center}
Abstract
\end{center}
{\itshape
We consider a viable market model. Suppose that new information arrives at the market. We are interested in modeling the market reaction facing to the change of information. In particular we seek for the limit on the intensity of information change below which the market stays always viable. We succeed to find such a limit with the drift operator when the market possesses the martingale representation property.
}

\

\section{Introduction}

There are situations where one is concerned with market models affected by the expansion of the information flow. It is the case, for example, when we consider the coherence between two parallel markets and we regard if the knowledge on one market brings arbitrage opportunities in the other one, or when we price defaultable corporate bonds and we are interested in the impact of the default information on the market of corporate stocks, or when we model the information asymmetry between different agents, or when we model the credit risk, etc. 

The fundamental tool to study such situations is the theory of enlargement of filtrations (see, for example, \cite{J, JY, Pr, mansuyYor, CJY}). Many works exist on market modeling in application of this theory. We observe, however, that the applications are essentially confined within two specific frameworks : the initial enlargement of filtrations or the progressive enlargement of filtrations, which fix the way that the information flow is expanded. In this paper we consider the problem from a different perspective. We will suppose the viability of the market after the expansion of the information flow. We regard then the consequences that we can draw about the expanded information flow. 

The viability of a market is the property that the utility optimization in the market have solutions (see \cite{ing, kabanov, loew, KC2010, K2012}). According to \cite{KC2010}, an operational way to formulate the market viability is to say that the price processes are semimartingales and they satisfy the condition \texttt{NA1}. The condition \texttt{NA1} is an operational condition, because it is equivalent to the existence of local martingale deflators and is equivalent to the structure condition. These last suit perfectly to stochastic calculus. (See for example \cite{Sch1, choulli, Sch2,K2012, Kar, takaoka, imkeller, Fon} for background information. Section \ref{deflator-al} gives precise definitions of these notions.)

An information flow in a market is modeled by a filtration $\mathbb{F}$. An expansion of the information flow is represented by another filtration $\mathbb{G}$ such that $\mathbb{F}\subset \mathbb{G}$. We are in the framework of enlargement of filtrations. A basic assumption in the theory of enlargement of filtrations is the Hypothesis$(H')$ (see \cite{J} or Section \ref{vocabulary}), which states that any real valued local martingale $X$ in $\mathbb{F}$ is a semimartingale in $\mathbb{G}$. Let $\Gamma(X)$ denote the drift part of $X$ in $\mathbb{G}$.  

The concept of information is a fascinating notion, but also an elusive notion especially when we want to quantify it. The framework of enlargement of filtrations $\mathbb{F}\subset\mathbb{G}$ offers since long a nice laboratory to test the ideas. In general, no common consensus exists how to quantify the difference between the two filtrations $\mathbb{G}$ and $\mathbb{F}$. The notion of entropy has been used there (see for example \cite{Y, ankirchnerImkeller}). But a more convincing measurement of information should be simply the drift operator $\Gamma$ itself. (See, for example, the discussion in \cite{ankirchner} and the references therein. See also \cite{JS} for a use of $\Gamma$ in a study of the martingale representation property in $\mathbb{G}$.) This observation is strengthened by the result in this paper. For the problem we study below, the drift operator $\Gamma(X)$ proves to be a very good gauge to control the level of the information expansion, in order to maintain the viability of the market. More precisely, suppose that the market with information flow $\mathbb{F}$ is viable on a time horizon $[0,T]$. Suppose that the Hypothesis$(H')$ holds for the passage from $\mathbb{F}$ into $\mathbb{G}$ with a drift operator $\Gamma(X)$ in the classical form :$$
\Gamma(X)=\transp\overline{\varphi}\centerdot[N,X]^{\mathbb{F}-p}
$$
where $N$ is a vector valued $\mathbb{F}$ local martingale, $\overline{\varphi}$ is a vector valued $\mathbb{G}$ predictable process, and $\bullet^{\mathbb{F}-p}$ denotes the predictable dual projection in $\mathbb{F}$. Then, under technical assumptions, if the increasing process (defined by an integral) $
\transp\overline{\varphi}(\centerdot[N^c,\transp N^c])\overline{\varphi}
$
is a finite process, if $1+\transp\overline{\varphi}\Delta N\geq \mathsf{u}$ for a $\mathbb{G}$ predictable process $\mathsf{u}$ and the increasing process 
$
\frac{1}{\mathsf{u}}\centerdot[D^d,D^d]$ and $\frac{1}{\mathsf{u}}\transp\overline{\varphi}(\centerdot[N^d,\transp N^d])\overline{\varphi}
$
are $(\mathbb{P},\mathbb{G})$ locally integrable, the market with the expanded information flow $\mathbb{G}$ will be viable on $[0,T]$. See Theorem \ref{main} and Corollary \ref{main-consequence} for the exact statements of the results.

If no jumps occurs in $\mathbb{F}$, the problem of the viability in the market with expanded information flow will have a very quick solution. The situation becomes radically different when jumps occur. To have an idea about the implication of jumps in the study of the market viability, the articles \cite{KC2010, K2012} give a good illustration. Technically speaking, this work is a study on the jumps of filtrations. We come to a satisfactory solution for the problem of jumps, under two key assumptions : The drift operator $\Gamma$ takes a classical form and the filtration $\mathbb{F}$ possesses the martingale predictable representation property. See Assumption \ref{assump1} and \ref{assump-mrt} for details.
 
The paper is organized as follows : In section \ref{vocabulary} we recall some vocabulary from stochastic calculus used in this paper. In section \ref{SC-def} we introduce the notion of structure condition and the basic assumptions. We give a first analysis on the problem to be solved in dividing the structure condition into three types : continuous type, accessible type and totally inaccessible type. We solve the problem in the case of continuous structure condition. In section \ref{sc+mrt}, noting that the two structure conditions of jumping type (accessible or totally inaccessible types) need a special treatment, we introduce the martingale representation property and we give a different formulation of the problems to be solved. The main result of this paper is stated in Theorem \ref{main}. The long section \ref{proof} is devoted to the proof of the main result. We note firstly that the structure condition of jumping type can be localized at the jumping times. This observation leads us to the following idea : solve the problems locally at each jumping times, and integrate the local solution into a global one. To find solutions at the jumping times, we need to compare precisely the conditional expectations $\mathbb{E}[\cdot|\mathcal{F}_{R-}]$ and $\mathbb{E}[\cdot|\mathcal{G}_{R-}]$ for a $\mathbb{F}$ stopping time $R$. This is done with the notion of conditional multiplicity introduced in \cite{BEKSY}. We obtain in subsection \ref{exp-at-jump} a nice relationship between these conditional expectations, which provide the key element to construct the final solutions.

We say two words on the role played by the martingale representation property (\mrt) in this paper. It is a tool for computations. This observation becomes particularly clear when we note that the essential condition in the main result Theorem \ref{main} is the integrability condition concerning the processes $\mathsf{u}, \overline{\varphi}, D, N$, which can be formulated independently without \mrt. This may mean two things : Firstly, we have the freedom to choose different \mrt\ to facilitate the computations. Secondly, the result of this paper may be true without the martingale representation property. For the moment, we are content to work with \mrt, because it leads to a complete solution with elementary computations. This founds a good basis for a future study.

\

\section{Vocabulary from stochastic calculus}\label{vocabulary}

This paper is based on the stochastic calculus as presented in \cite{Jacodlivre, JacShi, HWY}. We fix a probability space $(\Omega, \mathcal{A},\mathbb{P})$. Let $\mathbb{F}=(\mathcal{F}_t)_{t\geq 0}$ be a filtration in $\mathcal{A}$ satisfying the usual condition.

\

\textbf{Vectorial convention}

The elements $v$ in $\mathbb{R}^k$ are considered as vertical vectors. We denote their transpositions by $\transp v$. We denote by $(\!(v|v)\!)$ the inner product in $\mathbb{R}^k$. 

We deal with finite family of real processes $X=(X_i)_{1\leq i\leq k}$ ($k\in\mathbb{N}^*$). They will be considered as a process $X$ taking values in the vector space $\mathbb{R}^k$. To mention such a $X$, we say that $X$ is a $k$-dimensional process. For the value of $X$ at time $t\geq 0$, we denote by $(X_t)_{i}$ the $i^e$ component of the vector $X_t$. When $X$ is a semimartingale, we denote by $[X,\transp X]$ the $k\times k$-dimensional matrix valued process whose coefficients are $[X_i,X_j]$ for $1\leq i,j\leq k$.

\

\textbf{The processes}

By definition, $\Delta_0X=0$ for any càdlàg process $X$. For an process $A$ with finite variation (always assumed càdlàg), we denote by $\mathsf{d}A$ the (signed) random measure that $A$ generates. For $p\geq 1$, for a $k$-dimensional process $V$ whose components are all processes with finite variation, we introduce$$
\dcb
\underline{\mathsf{V}}^p(\mathbb{P},\mathbb{F}, \mathsf{d}V)
=\{\overline{H} &:& \mbox{$\overline{H}$ is $k$-dimensional $\mathbb{F}$ predictable process, and for $1\leq h\leq k$,}\\ &&\mbox{  $\int_0^t|(\overline{H}_s)_h| |\mathsf{d}(V_s)_h|, t\geq 0,$ is a $(\mathbb{P},\mathbb{F})$ locally $p$-integrable process}\}
\dce
$$
When $p=1$, we note simply $\underline{\mathsf{V}}(\mathbb{P},\mathbb{F}, \mathsf{d}V)$.

\pagebreak

\textbf{The projections}

With respect to the filtration $\mathbb{F}$, the notation ${^{\mathbb{F}-p}}\bullet$ denotes the predictable projection, and the notation $\bullet^{\mathbb{F}-p}$ denotes the predictable dual projection.

\textbf{The martingales and the semimartingales}

For any special semimartingale $X$, we can decompose $X$ in the form (see \cite[Theorem 7.25]{HWY}) :$$
\dcb
X=X_0+X^m+X^v\\
X^m=X^c+X^{da}+X^{di}\\
X^v=X^{vc}+X^{vj}
\dce
$$
where $X^m$ is the martingale part of $X$ and $X^v$ is the predictable part of finite variation of $X$, $X^c$ is the continuous martingale part, $X^{da}$ is the part of compensated sum of accessible jumps, $X^{di}$ is the part of compensated sum of totally inaccessible jumps, $X^{vc}$ is the continuous part of $X^v$, $X^{vj}$ is the purely jumping part of $X^v$. We recall that this decomposition of $X$ depends on the reference probability and the reference filtration. In the computations below we apply this notation system only for the decompositions in $\mathbb{F}$. We recall that every part of the decomposition of $X$, except $X_0$, is assumed null at $t=0$. 

For a semimartingale $X$, $[X,X]_0=0$ by definition. Let $p\geq 1$, we denote by $\mathcal{H}^p(\mathbb{P},\mathbb{F})$ the space of $(\mathbb{P},\mathbb{F})$ local martingales $X$ such that $\sqrt{[X,X]_\infty}^p$ is $(\mathbb{P},\mathbb{F})$ integrable. We denote by $\mathcal{H}^p_{loc}(\mathbb{P},\mathbb{F})$ the family of local martingales locally in $\mathcal{H}^p(\mathbb{P},\mathbb{F})$. We say the a sequence $(X_n)$ in the space $\mathcal{H}^p_{loc}(\mathbb{P},\mathbb{F})$ converges to $X$, if there exists an increasing sequence $(R_m)$ of $\mathbb{F}$ stopping times such that $\sup_{m\geq 1}R_m=\infty$ and, for each $m$, $X_n^{R_m}$ converges to $X^{R_m}$ in $\mathcal{H}^p(\mathbb{P},\mathbb{F})$.

\

\textbf{The stochastic integrals}

In this paper we employ the notion of stochastic integral only on the predictable processes. The stochastic integral are defined as 0 at $t=0$. We use a point "$\centerdot$" to indicate the integrator process in a stochastic integral. For example, the stochastic integral of a real predictable process $\overline{H}$ with respect to a real semimartingale $Y$ is denoted by $\overline{H}\centerdot Y$, while the expression $\transp\overline{K}(\centerdot[X,\transp X])\overline{K}$ denotes the process$$
\int_0^t \sum_{i=1}^k\sum_{j=1}^k(\overline{K}_s)_{i,s}(\overline{K}_s)_{j,s} d[X_i,X_j]_s,\ t\geq 0,
$$
where $\overline{K}$ is a $k$-dimensional predictable process and $X$ is a $k$-dimensional semimartingale. The expression $\transp\overline{K}(\centerdot[X,\transp X])\overline{K}$ respects the matrix product rule. The value at $t\geq 0$ of a stochastic integral will be denoted, for example, by $\transp\overline{K}(\centerdot[X,\transp X])\overline{K}_t$.

The notion of the stochastic integral with respect to a $k$-dimensional local martingale $X$ follows \cite{Jacodlivre}. We introduce $\underline{\mathsf{L}}^p(\mathbb{P},\mathbb{F},X)$ as the family of $\mathbb{F}$ predictable $k$-dimensional process $\overline{H}$ such that $\sqrt{\transp\overline{H}(\centerdot[X,\transp X])\overline{H}}$ is $(\mathbb{P},\mathbb{F})$ locally $p$-times integrable $(p\geq 1)$. For a sequence $(\overline{H}_n)_{n\geq 1}$ of elements in $\underline{\mathsf{L}}^p(\mathbb{P},\mathbb{F},X)$, we say that $\overline{H}_n$ converges to an element $\overline{H}\in\underline{\mathsf{L}}^p(\mathbb{P},\mathbb{F},X)$, if there exists an increasing sequence $(T_m)_{m\geq 1}$ of $\mathbb{F}$ stopping times such that $\sup_{m\geq 1}T_m=\infty$ and, for any $m\geq 1$, 
$$
\mathbb{E}_{\mathbb{P}}[\sqrt{\transp(\overline{H}-\overline{H}_n)(\centerdot[X,\transp X])(\overline{H}-\overline{H}_n)}^p_{T_m}]\longrightarrow 0,\ n\uparrow\infty.
$$
We define in the same way the notion of Cauchy sequence in $\underline{\mathsf{L}}^p(\mathbb{P},\mathbb{F},X)$. We recall that the space $\underline{\mathsf{L}}^p(\mathbb{P},\mathbb{F},X)$ is complete in the sens that for any Cauchy sequence $(\overline{H}_n)_{n\geq 1}$ in $\underline{\mathsf{L}}^p(\mathbb{P},\mathbb{F},X)$, there exists a $\overline{H}\in\underline{\mathsf{L}}^p(\mathbb{P},\mathbb{F},X)$ such that $\overline{H}_n$ converges to $\overline{H}$ (see \cite[Théorème (4.60)]{Jacodlivre}). When $p=1$, we denote $\underline{\mathsf{L}}(\mathbb{P},\mathbb{F},X)$ instead of $\underline{\mathsf{L}}^1(\mathbb{P},\mathbb{F},X)$. We say that a $k$-dimensional $\mathbb{F}$ predictable process is integrable with respect to $X$ under the probability $\mathbb{P}$, if it belongs to $\underline{\mathsf{L}}(\mathbb{P},\mathbb{F},X)$. For such an integrable process $\overline{H}$, the stochastic integral $\transp\overline{H}\centerdot X$ is defined. We introduce $\mathfrak{L}(\mathbb{P},\mathbb{F},X)$ the set of all stochastic integrals $\transp\overline{H}\centerdot X$ for $\overline{H}\in\underline{\mathsf{L}}(\mathbb{P},\mathbb{F},X)$. The bracket process of a stochastic integral $\transp\overline{H}\centerdot X$ can be computed using Remarque(4.36) and Proposition(4.68) in \cite{Jacodlivre}. We note that, if $\overline{H}, \overline{K}\in\underline{\mathsf{L}}^2(\mathbb{P},\mathbb{F},X)$, the process $\transp \overline{H}(\centerdot[X,\transp X])\overline{K}$ is a well defined process with $(\mathbb{P},\mathbb{F})$ locally integrable variation.

\

\textbf{The martingale representation property}

We introduce the space $\mathfrak{L}(\mathbb{P},\mathbb{F})$ of all $(\mathbb{P},\mathbb{F})$ local martingales $M$ null at the origin $M_0=0$. We consider a $k$-dimensional stochastic process $W$. We say that the martingale representation property holds in the filtration $\mathbb{F}$ under the probability $\mathbb{P}$ with respect to the driving process $W$, if $W$ is a $(\mathbb{P},\mathbb{F})$ local martingale, and if $
\mathfrak{L}(\mathbb{P},\mathbb{F},W) = \mathfrak{L}(\mathbb{P},\mathbb{F}).
$
The martingale representation property will be denoted by \mrt$(\mathbb{P},\mathbb{F},W)$, or simply by \mrt. See \cite{JacShi} for general information on \mrt.

\

\textbf{Enlargements of filtrations and the Hypothesis$(H')$ }

Let $\mathbb{G}$ be a filtration containing $\mathbb{F}$. Let $T$ be a $\mathbb{G}$ stopping time. We introduce the Hypothesis$(H')$:   

\textit{
\noindent\textbf{Hypothesis$(H')$ on the time horizon $[0,T]$} There exists an application $\Gamma$ from $\mathfrak{L}(\mathbb{P},\mathbb{F})$ into the space of càdlàg $\mathbb{G}$-predictable processes on $[0,T]$ with finite variation and null at the origin, such that, for any $X\in \mathfrak{L}(\mathbb{P},\mathbb{F})$, $\widetilde{X}:=X-\Gamma(X)$ is a $(\mathbb{P},\mathbb{G})$ local martingale on $[0,T]$. The operator $\Gamma$ will be called the drift operator.
}

\textbf{Other conventions}

In this paper, calling a number $a$ positive means that $a\geq 0$ and calling a function $f(t),t\in\mathbb{R},$ an increasing function means $f(s)\leq f(t)$ for $s\leq t$. For a number $a>0$, we call it a strictly positive number.

Relations between random variables is to be understood almost sure relations. For a random variable $X$ and a $\sigma$-algebra $\mathcal{F}$, the expression $X\in\mathcal{F}$ means that $X$ is $\mathcal{F}$-measurable. 

\

\section{Structure condition}\label{SC-def}

\subsection{Basic setting and definitions}\label{deflator-al}

All long this paper, $S=(S_1,\ldots, S_\mathsf{k})$ denotes a $\mathsf{k}$-dimensional $(\mathbb{P},\mathbb{F})$ special semimartingale, whose components are all strictly positive. The triple $(\mathbb{P},\mathbb{F},S)$ represents a financial market with a discounted price process $S$ and an information flow $\mathbb{F}$. Let $x>0$. We introduce $$
\dcb
\mathfrak{A}_x(\mathbb{P},\mathbb{F},S)=\{\overline{H}&:& \overline{H}\in\underline{\mathsf{L}}(\mathbb{P},\mathbb{F},S^m)\cap\underline{\mathsf{V}}(\mathbb{P},\mathbb{F},\mathsf{d}S^v),\ (x+\transp\overline{H}\centerdot S)_t \geq 0, \forall t\geq 0\}
\dce
$$
A couple $\pi=(x,\overline{H})$, where $x>0$ and $\overline{H}\in\mathfrak{A}_x(\mathbb{P},\mathbb{F},S)$, will be called an admissible strategy in the market $(\mathbb{P},\mathbb{F},S)$. We set $
S^{\pi}=S^{(x,\overline{H})}=x+\transp\overline{H}\centerdot S.
$

\bd
Let $T>0$ be a $\mathbb{F}$ stopping time. We say that the price process $S$ satisfies the structure condition in the filtration $\mathbb{F}$ with a $(\mathbb{P},\mathbb{F})$ local martingale $D$ under the probability $\mathbb{P}$ on the time horizon $[0,T]$, if there exists a real $(\mathbb{P},\mathbb{F})$ local martingale $D$ such that, on the time interval $[0,T]$, $D_0=0, \Delta D<1$, $[S^m_{i}, D]^{\mathbb{F}-p}$ exists, and $S^{v}_{i} = [S^m_{i}, D]^{\mathbb{F}-p}$ for $1\leq i\leq \mathsf{k}$. 
\ed

\bd
Let $T>0$ be a $\mathbb{F}$ stopping time. We call a strictly positive real process $\Upsilon$ with $\Upsilon_0=1$ a local martingale deflator on the time horizon $[0,T]$ for the market $(\mathbb{P},\mathbb{F},S)$, if, for any admissible strategy $(x,\overline{H})$, the process $\Upsilon S^{(x,\overline{H})}$ is a $(\mathbb{P},\mathbb{F})$ local martingale on $[0,T]$. In particular, $\Upsilon$ is a $(\mathbb{P},\mathbb{F})$ local martingale. 
\ed

\brem
Note that the notion of structure condition exists in the literature, particularly in \cite{choulli, Sch2}. In this paper we have formulated the structure condition in a form slightly different from the original one, in order to better adapt to the problem studied below.
\erem

The existence of local martingale deflator and the structure condition are equivalent conditions, as it will be explained below. They are equivalent also to the conditions \texttt{NUPBR} and \texttt{NA1} (see \cite{KC2010, takaoka}). 

\bethe
Let $T>0$ be a $\mathbb{F}$ stopping time. The market $(\mathbb{P},\mathbb{F},S)$ possesses a local martingale deflator on the time horizon $[0,T]$, if and only if $S$ satisfies the structure condition on the time horizon $[0,T]$.
\ethe

\textbf{Proof.} It is a standard application of the integration by parts formula, once we know that a local martingale deflator is always an exponential local martingale, which can be proved using \cite{Jacodlivre} or \cite{choulli2}. \ok

We assume in the remainder part of this paper the following assumption.

\bassumption\label{assump0}
$S$ satisfies the structure condition in $(\mathbb{P},\mathbb{F})$ on the time horizon $[0,\infty)$ with a $(\mathbb{P},\mathbb{F})$ local martingale $D$. 
\eassumption

\

\subsection{Structure condition in a market with an expanded information flow} 

We look at the situation that new information arrives at the market $(\mathbb{P},\mathbb{F},S)$. This situation is modeled by $(\mathbb{P},\mathbb{G},S)$ with $\mathbb{G}$ being an expansion of the filtration $\mathbb{F}$ : $\mathcal{G}_t\supset\mathcal{F}_t$ for $t\geq 0$. We now study the structure condition in this expanded market $(\mathbb{P},\mathbb{G},S)$. 

We assume the following conditions :

\bassumption\label{assump1} Let $T$ be a $\mathbb{G}$ stopping time.
\
\ebe 
\item
The Hypothesis$(H')$ is satisfied on the time horizon $[0,T]$ with a drift operator $\Gamma$.

\item
There exist $N=(N_1,\ldots,N_n)$ an $n$-dimensional $(\mathbb{P},\mathbb{F})$ local martingale, and $\overline{\varphi}$ an $n$ dimensional $\mathbb{G}$ predictable process such that, for any $X\in\mathfrak{L}(\mathbb{P},\mathbb{F})$, $[N,X]^{\mathbb{F}-p}$ exists, $\overline{\varphi}\in\underline{\mathsf{V}}(\mathbb{P},\mathbb{G},\mathsf{d}[N,X]^{\mathbb{F}-p})$, and $$
\Gamma(X)=\transp\overline{\varphi}\centerdot [N,X]^{\mathbb{F}-p}
$$
on the time horizon $[0,T]$. 
\dbe
\eassumption

\brem
Assumption \ref{assump1} is satisfied in a number of concrete models, especially among the models of progressive enlargement of filtrations, or the models satisfying Jacod's criterion (\cite{Jacod}), or $(\natural)$-model in \cite{JS}. Assumption \ref{assump1} is a fairly strong condition. It assumes almost that $N$ is locally a \texttt{BMO} martingale ($\overline{\varphi}$ may contain singularity). The purpose of this paper is to understand the structure condition in the expanded market at least under this strong condition. \ok 
\erem

Note that, even if the Hypothesis$(H')$ holds only on the time horizon $[0,T]$, it is convenient to extend the operator $\Gamma$ on the whole $\mathbb{R}_+$. For this, we simply suppose that $\overline{\varphi}=\overline{\varphi}\ind_{[0,T]}$ and define $\Gamma$ by $\transp\overline{\varphi}\centerdot [N,X]^{\mathbb{F}-p}$ (recall that the processes $N, \overline{\varphi}$ are defined on the whole $\mathbb{R}_+$). As a consequence, $\widetilde{X}=X-\Gamma(X)$ and $[\widetilde{X},\widetilde{X}]$ for $X\in\mathfrak{L}(\mathbb{P},\mathbb{F})$ are well defined on $\mathbb{R}_+$.

Under Assumptions \ref{assump0} and \ref{assump1}, we have immediately the following lemma.

\bl
We denote $M=S^m$. The $(\mathbb{P},\mathbb{G})$ canonical decomposition of $S$ on $[0,T]$ is given by $$
S=\widetilde{M}+[D,M]^{\mathbb{F}-p}+\transp\overline{\varphi}\centerdot [N,M]^{\mathbb{F}-p}
$$
so that the structure condition for the expanded market $(\mathbb{P},\mathbb{G},S)$ takes the following form : There exists a $\mathbb{G}$ local martingale $Y$ such that $Y_0=0, \Delta Y<1$, $[Y,\widetilde{M}_i]^{\mathbb{G}-p}$ exists, and for $1\leq i\leq \mathsf{k}$, 
\begin{equation}\label{structure-condition}
[Y,\widetilde{M}_i]^{\mathbb{G}-p}=[D,M_i]^{\mathbb{F}-p}+\transp\overline{\varphi}\centerdot [N,M_i]^{\mathbb{F}-p}
\end{equation}
on the time horizon $[0,T]$.
\el

\

\subsection{Structure condition decomposed}

The structure condition (\ref{structure-condition}) is naturally linked with the following more specific structure conditions :  
\ebe
\item[. ]\textbf{Continuous structure condition}
There exists a continuous $(\mathbb{P},\mathbb{G})$ local martingale $Y'$ such that, on the time horizon $[0,T]$,
\begin{equation}\label{structure-condition-c}
[Y',\widetilde{M^{c}_i}]^{\mathbb{G}-p}
=[D,M^{c}_i]^{\mathbb{F}-p}+\transp\overline{\varphi}\centerdot [N,M^{c}_i]^{\mathbb{F}-p}
\end{equation}
\vspace{13pt}

\item[. ]\textbf{Accessible structure condition}
There exists a $(\mathbb{P},\mathbb{G})$ local martingale $Y''$ of compensated sum of jumps such that $[Y'',\widetilde{M^{da}_i}]^{\mathbb{G}-p}$ exists; the jumps of $Y''$ are all $(\mathbb{P},\mathbb{F})$ accessible; $\Delta Y''<1$; and, on the time horizon $[0,T]$,
\begin{equation}\label{structure-condition-da}
[Y'',\widetilde{M^{da}_i}]^{\mathbb{G}-p}
=[D,M^{da}_i]^{\mathbb{F}-p}+\transp\overline{\varphi}\centerdot [N,M^{da}_i]^{\mathbb{F}-p}
\end{equation}
\vspace{13pt}

\item[. ]\textbf{Totally inaccessible structure condition}
There exists a $(\mathbb{P},\mathbb{G})$ local martingale $Y'''$ of compensated sum of jumps such that $[Y''',\widetilde{M^{di}_i}]^{\mathbb{G}-p}$ exists; the jumps of $Y'''$ are all $(\mathbb{P},\mathbb{F})$ totally inaccessible; $\Delta Y''<1$; and, on the time horizon $[0,T]$,
\begin{equation}\label{structure-condition-di}
[Y''',\widetilde{M^{di}_i}]^{\mathbb{G}-p}
=[D,M^{di}_i]^{\mathbb{F}-p}+\transp\overline{\varphi}\centerdot [N,M^{di}_i]^{\mathbb{F}-p}
\end{equation}
\vspace{17pt}
\dbe

We recall that Assumptions \ref{assump0} and \ref{assump1} are in force.

\bl
The structure condition (\ref{structure-condition}) implies the condition (\ref{structure-condition-da}).
\el

\textbf{Proof.} The lemma follows from the observation that there exists a $\mathbb{F}$-predictable real process $K$ such that $$
M^{da}_i=K\centerdot M_i \mbox{ and } \widetilde{M^{da}_i}=K\centerdot \widetilde{M_i}, \ 1\leq i\leq \mathsf{k}\ \ok
$$

\bl\label{2341}
The structure conditions (\ref{structure-condition-c}) and (\ref{structure-condition-da}) and (\ref{structure-condition-di}) imply the structure condition (\ref{structure-condition}). 
\el

\textbf{Proof.} Let $Y', Y'', Y'''$ be respectively the solutions of the structure conditions (\ref{structure-condition-c}) and (\ref{structure-condition-da}) and (\ref{structure-condition-di}). There exist strong orthogonalities between $Y', Y'', Y'''$ and $\widetilde{M^{c}},\widetilde{M^{da}}, \widetilde{M^{di}}$:$$
[Y' +Y'' +Y''', \widetilde{M^{c}}+\widetilde{M^{da}}+ \widetilde{M^{di}}]
=[Y', \widetilde{M^{c}}]+[Y'', \widetilde{M^{da}}]
+[Y''',  \widetilde{M^{di}}]
$$
With these orthogonalities, we check immediately that $Y'+Y''+Y'''$ is a solution of the structure condition (\ref{structure-condition}). \ok 

\

\bethe\label{answer-c}
We denote by $\overline{p}_D\centerdot M^c$ and $\overline{p}_N\centerdot M^c$ the $(\mathbb{P},\mathbb{F})$ orthogonal projection of $D^c$ and respectively of $N^c$ on the stable space generated by the $(\mathbb{P},\mathbb{F})$ local martingales $M^c_i, 1\leq i\leq \mathsf{k}$, where $\overline{p}_D$ is a vector valued process in $\underline{\mathsf{L}}(\mathbb{P},\mathbb{F},M^c)$ and $\overline{p}_N$ is a matrix valued process whose line vectors are in $\underline{\mathsf{L}}(\mathbb{P},\mathbb{F},M^c)$.  
Then, the structure condition (\ref{structure-condition-c}) is satisfied, 
if and only if $
\transp\overline{\varphi}(\centerdot[N^c,\transp N^c])\overline{\varphi}
$
is a finite process. In this case, $Y'=(\overline{p}_D+\transp\overline{\varphi}\overline{p}_N)\centerdot \widetilde{M^c}^T$ satisfies the structure condition (\ref{structure-condition-c}).
\ethe

\textbf{Proof.} For the existence of $\overline{p}_D$ and $\overline{p}_N$, see Théorème (4.35) and comment on Proposition (4.26) in \cite{Jacod}. 

Let us consider the following equation on $[0,T]$:
$$
\transp \overline{K}\centerdot [\widetilde{M^{c}},\transp\widetilde{M^{c}}]^{\mathbb{G}-p}
=\transp\overline{p}_D\centerdot [M^{c},\transp M^{c}]^{\mathbb{F}-p}
+\transp\overline{\varphi}\overline{p}_N\centerdot [M^{c},\transp M^{c}]^{\mathbb{F}-p}
$$
for a $\mathbb{G}$ predictable process $\overline{K}$. We note that, on $[0,T]$,$$
[\widetilde{M^{c}},\transp\widetilde{M^{c}}]^{\mathbb{G}-p}
=[M^{c},\transp M^{c}]^{\mathbb{G}-p}
=[M^{c},\transp M^{c}]^{\mathbb{F}-p}
$$
The equation becomes 
$$
\transp \overline{K}\centerdot [\widetilde{M^{c}},\transp\widetilde{M^{c}}]^{\mathbb{G}-p}
=(\transp\overline{p}_D
+\transp\overline{\varphi}\overline{p}_N)\centerdot [\widetilde{M^{c}},\transp\widetilde{M^{c}}]^{\mathbb{G}-p}
$$
The solution of this equation can only be $$
K=\overline{p}_D+\transp\overline{\varphi}\overline{p}_N,\ \  \mathsf{d}[\widetilde{M^{c}},\transp\widetilde{M^{c}}]^{\mathbb{G}-p} - a.s..
$$
on $[0,T]$. 

Now, if the structure condition (\ref{structure-condition-c}) holds, by orthogonal projection, we can replace $Y'$ in  equation (\ref{structure-condition-c}) by a local martingale of the form $\transp\overline{K}\centerdot \widetilde{M^c}^T$ with $\overline{K}\in\underline{\mathsf{L}}(\mathbb{P},\mathbb{G},\widetilde{M^c}^T)$. We see that this $\overline{K}$ has to be a solution of the above equation. Consequently, $\overline{p}_D+\transp\overline{\varphi}\overline{p}_N=\overline{K}\in\underline{\mathsf{L}}(\mathbb{P},\mathbb{G},\widetilde{M^c}^T)$. This is equivalent to say that $\transp\overline{\varphi}(\centerdot[N^c,\transp N^c])\overline{\varphi}$ is a finite process.

On the other hand, if $\overline{p}_D+\transp\overline{\varphi}\overline{p}_N\in\underline{\mathsf{L}}(\mathbb{P},\mathbb{G},\widetilde{M^c}^T)$, the local martingale $Y'=(\overline{p}_D+\transp\overline{\varphi}\overline{p}_N)\centerdot \widetilde{M^c}^T$ will satisfy the structure condition (\ref{structure-condition-c}). The theorem is proved. \ok

\

\section{Structure conditions (\ref{structure-condition-da}) and (\ref{structure-condition-di}) under the martingale representation property}\label{sc+mrt}

For the question when a structure condition holds, Theorem \ref{answer-c} gives a complete answer in the case of structure conditions (\ref{structure-condition-c}) with a simple proof thanks to the continuity. In the contrast, the jumping nature of the conditions (\ref{structure-condition-da}) and (\ref{structure-condition-di}) makes the question more cumbersome, notably because of the difference between the $\mathbb{G}$ predictable bracket and the $\mathbb{F}$ predictable bracket $[\widetilde{M},\transp \widetilde{M}]^{\mathbb{G}-p}$ and $[M,\transp M]^{\mathbb{F}-p}$. 

In this section, using the martingale representation property, we rewrite the structure conditions (\ref{structure-condition-da}) and (\ref{structure-condition-di}) in a different form which will be used later to find their solutions.  

\subsection{Computing predictable dual projections in $\mathbb{G}$}

We recall Assumptions \ref{assump0} and \ref{assump1} being in force. As a consequence, we have immediately the following lemma.

\bl\label{A-G-p}
For any $\mathbb{F}$ adapted càdlàg process $A$ with $(\mathbb{P},\mathbb{F})$ locally integrable variation, we have $$
A^{\mathbb{G}-p}=A^{\mathbb{F}-p}+\Gamma(A-A^{\mathbb{F}-p})=A^{\mathbb{F}-p}+\transp\overline{\varphi}\centerdot[N,A-A^{\mathbb{F}-p}]^{\mathbb{F}-p}
$$
on $[0,T]$. In particular, for $R$ a $\mathbb{F}$ stopping time either $(\mathbb{P},\mathbb{F})$ predictable or $(\mathbb{P},\mathbb{F})$ totally inaccessible, for $\xi\in\mathbb{L}^1(\mathbb{P},\mathcal{F}_{R})$, $$
(\xi\ind_{[R,\infty)})^{\mathbb{G}-p}
=
(\xi\ind_{[R,\infty)})^{\mathbb{F}-p}+\transp\overline{\varphi}\centerdot(\Delta_{R}N\xi\ind_{[R,\infty)})^{\mathbb{F}-p}
$$
on $[0,T]$.
\el

\

\subsection{\mrt\ assumption}

We suppose from now on another assumption. 

\bassumption\label{assump-mrt}
\
\ebe
\item
$M=S^{m}$, $D$ and $N$ are $(\mathbb{P},\mathbb{F})$ locally square integrable local martingales.  
\item
\mrt$(\mathbb{P},\mathbb{F},W)$ where $W$ is a $d$-dimensional $(\mathbb{P},\mathbb{F})$ locally square integrable local martingale. 
\dbe
\eassumption

Let $\overline{K}$ be a $d$-dimensional vector valued $\mathbb{G}$ predictable process in $\underline{\mathsf{L}}^2(\mathbb{P},\mathbb{G},\widetilde{W^{da}}^T)$ (resp. in $\underline{\mathsf{L}}^2(\mathbb{P},\mathbb{G},\widetilde{W^{di}}^T)$). We will say that $\overline{K}$ solves the equation (\ref{structure-condition-da}), (resp. the equation (\ref{structure-condition-di})) if $Y=\overline{K}\centerdot \widetilde{W^{da}}^T$ (resp. $Y=\overline{K}\centerdot \widetilde{W^{di}}^T$) satisfies the equation (\ref{structure-condition-da}) (resp. the equation (\ref{structure-condition-di}))

Under Assumption \ref{assump-mrt}, any $(\mathbb{P},\mathbb{F})$ local martingale $X$ is of the form $\transp\overline{k}\centerdot W$, where $\overline{k}\in\underline{\mathsf{L}}(\mathbb{P},\mathbb{F},W)$. We call the process $\overline{k}$ the coefficient of $X$ in its martingale representation with respect to the driving process $W$. This appellation extends naturally to vector valued local martingales. There exists in general no uniqueness for the coefficients.

Let be respectively the coefficients of $M$, of $D$ and of $N$, the $\mathsf{k}\times d$ dimensional matrix valued process $\overline{m}$, the $d$-dimensional vector valued process $\overline{d}$, and the $n\times d$-dimensional matrix valued process $\overline{\zeta}$. 

\brem
A \mrt\ condition is a fairly strong condition. However, the basic idea of this paper is to investigate the impact with which new information can affect an existing good market. \mrt\ is a reasonable condition on a good market. \ok
\erem

\subsection{Equation (\ref{structure-condition-da}) detailed through \mrt}

We are going to transcribe the equations (\ref{structure-condition-da}) and (\ref{structure-condition-di}) in term of the coefficients $\overline{m}, \overline{d}$ and $\overline{\zeta}$. We define a $d\times d\times d$ dimensional $\mathbb{F}$-predictable process $\overline{\kappa}$ by the following relation through the \mrt\ property:
$$
[W^{da}_a,W^{da}_b]-[W^{da}_a,W^{da}_b]^{\mathbb{F}-p}
=\sum_{e=1}^d\overline{\kappa}_{a,b,e}\centerdot W_e
$$
i.e. for $1\leq a,b\leq d$, the vector valued process $\overline{\kappa}_{a,b,\cdot}$ is the coefficient of $[W^{da}_a,W^{da}_b]-[W^{da}_a,W^{da}_b]^{\mathbb{F}-p}$ with respect to $W$. Note that necessarily$$
\transp\overline{\kappa}_{a,b,\cdot}\centerdot W
=\transp\overline{\kappa}_{a,b,\cdot}\centerdot W^{da}.
$$
Let $a$ be the matrix valued process defined by $a_{j,i}=[W^{da}_j,W^{da}_i]^{\mathbb{F}-p}, 1\leq j,i\leq d$. We have obviously $\overline{\kappa}_{b,a,e}=\overline{\kappa}_{a,b,e}$ and $a_{b,a}=a_{a,b}$, but also$$
\dcb
\sum_{e=1}^d\overline{\kappa}_{a,b,e}\centerdot a_{e,c}
&=&\sum_{e=1}^d\overline{\kappa}_{c,a,e}\centerdot a_{e,b}\\
\dce
$$
Let $\overline{\mathsf{f}}=\transp\overline{\zeta}\overline{\varphi}$ i.e., the $d$-dimensional vector valued process whose $e^{\mbox{\scriptsize th}}$ component is given by $
\sum_{b=1}^n
\overline{\varphi}_{b}\overline{\zeta}_{b,e}
$.
Let $\overline{\mathsf{ce}}$ denote the $d\times d$ dimensional matrix valued process whose element at row $j$ and at column $e$ is given by $$
\sum_{a=1}^n\sum_{b=1}^d\overline{\varphi}_{a}\overline{\zeta}_{a,b}\overline{\kappa}_{b,j,e}
=\sum_{b=1}^d\overline{\mathsf{f}}_{b}\overline{\kappa}_{b,j,e}
$$

Let $\mathfrak{I}_d$ denote the $d\times d$-dimensional identity matrix.

\bl\label{sc-da}
Let $\overline{K}\in\underline{\mathsf{L}}^2(\mathbb{P},\mathbb{G},\widetilde{W^{da}}^T)$. Then, $\overline{K}$ solves the equation (\ref{structure-condition-da}) if and only if
$$
\transp\overline{K} \left(\mathfrak{I}_d+\overline{\mathsf{ce}} - \Delta a \overline{\mathsf{f}}\transp\overline{\mathsf{f}}\right)(\centerdot a) \transp\overline{m}
=(\transp\overline{d}+\transp\overline{\mathsf{f}}) (\centerdot a) \transp\overline{m}\ 
$$
on the time horizon $[0,T]$.  
\el

\textbf{Proof.} With the new notations, $\overline{K}$ solves the equation (\ref{structure-condition-da})
if and only if, on $[0,T]$,$$
\transp\overline{K} (\centerdot [\widetilde{W^{da}},\transp\widetilde{W^{da}}]^{\mathbb{G}-p})\transp\overline{m}
=(\transp\overline{d}+\transp\overline{\mathsf{f}}) (\centerdot a) \transp\overline{m}
$$
On the time horizon $[0,T]$, for $1\leq j,i\leq d$,
$$\dcb
[\widetilde{W^{da}}_j,\widetilde{W^{da}_i}]^{\mathbb{G}-p}

&=&[W^{da}_j,W^{da}_i]^{\mathbb{G}-p}-[\Gamma(W^{da}_j),\Gamma(W^{da}_i)]^{\mathbb{G}-p}
\dce
$$ 
For the part concerning $\Gamma$, we have $$
\dcb
[\Gamma(W^{da}_j),\Gamma(W^{da}_i)]^{\mathbb{G}-p}

=\sum_{s\leq \cdot}\transp\overline{\varphi}_s\Delta_s[N,W^{da}_j]^{\mathbb{F}-p}\transp\overline{\varphi}_s\Delta_s[N,W^{da}_i]^{\mathbb{F}-p}

&=&\left((\Delta a\overline{\mathsf{f}}\transp\overline{\mathsf{f}})\centerdot a\right)_{j,i}\\

\dce
$$
Applying Lemma \ref{A-G-p} we compute also
$$
\dcb
[W^{da}_j,W^{da}_i]^{\mathbb{G}-p}
=[W^{da}_j,W^{da}_i]^{\mathbb{F}-p}+\Gamma([W^{da}_j,W^{da}_i]-[W^{da}_j,W^{da}_i]^{\mathbb{F}-p})

&=&\left((\mathfrak{I}_d+\overline{\mathsf{ce}})\centerdot a\right)_{j,i}
\dce
$$
Putting these together we write on $[0,T]$:
\begin{equation}\label{jcef}
[\widetilde{W^{da}},\transp\widetilde{W^{da}}]^{\mathbb{G}-p}
=\left(\mathfrak{I}_d+\overline{\mathsf{ce}}\right)\centerdot a - \Delta a \overline{\mathsf{f}}\transp\overline{\mathsf{f}}\centerdot a
\end{equation}
This proves the lemma. \ok

\subsection{Equation (\ref{structure-condition-di}) detailed through \mrt}

We define $\overline{\varrho}_{j,i}$ by
$$
[W^{di}_a,W^{di}_b]-[W^{di}_a,W^{di}_b]^{\mathbb{F}-p}
=\sum_{e=1}^d\overline{\varrho}_{a,b,e}\centerdot W_e
$$
i.e. for $1\leq a,b\leq d$, the vector valued process $\overline{\rho}_{a,b,\cdot}$ is the coefficient of $[W^{di}_a,W^{di}_b]-[W^{da}_a,W^{da}_b]^{\mathbb{F}-p}$ with respect to $W$. Note that necessarily$$
\transp\overline{\rho}_{a,b,\cdot}\centerdot W
=\transp\overline{\rho}_{a,b,\cdot}\centerdot W^{di}
$$
Let $b$ be the matrix valued process defined by $b_{j,i}=[W^{di}_j,W^{di}_i]^{\mathbb{F}-p}, 1\leq j,i\leq d$. We have $\overline{\varrho}_{b,a,e}=\overline{\varrho}_{a,b,e}$, $b_{b,a}=b_{a,b}$, and$$
\dcb
\sum_{e=1}^d\overline{\varrho}_{a,b,e}\centerdot b_{e,c}
&=&\sum_{e=1}^d\overline{\varrho}_{c,a,e}\centerdot b_{e,b}\\
\dce
$$
Let $\overline{\mathsf{ci}}$ denote the $d\times d$ dimensional matrix valued process whose element at row $j$ and at column $e$ is given by $$
\sum_{a=1}^n\sum_{b=1}^d\overline{\varphi}_{a}\overline{\zeta}_{a,b}\overline{\varrho}_{b,j,e}
=\sum_{b=1}^d\overline{\mathsf{f}}_{b}\overline{\varrho}_{b,j,e}
$$

\bl\label{sc-di}
Let $\overline{K}\in\underline{\mathsf{L}}^2(\mathbb{P},\mathbb{G},\widetilde{W^{di}}^T)$. Then, $\overline{K}$ solves the equation (\ref{structure-condition-di}), if and only if
$$
\transp\overline{K} \left(\mathfrak{I}_d+\overline{\mathsf{ci}}\right) (\centerdot b)\transp\overline{m}
=(\transp\overline{d}+\transp\overline{\mathsf{f}}) (\centerdot b) \transp\overline{m}
$$
on the time horizon $[0,T]$
\el

\textbf{Proof.} $\overline{K}$ solves the equation (\ref{structure-condition-di})
if and only if, on $[0,T]$,$$
\transp\overline{K} (\centerdot[\widetilde{W^{di}},\transp\widetilde{W^{di}}]^{\mathbb{G}-p})\transp\overline{m}
=(\transp\overline{d}+\transp\overline{\mathsf{f}}) (\centerdot b) \transp\overline{m}\ 
$$
We compute now $[\widetilde{W^{di}},\transp\widetilde{W^{di}}]^{\mathbb{G}-p}$ on $[0,T]$. Note that $\Gamma(W^{di})$ is continuous on $[0,T]$. For $1\leq i,j\leq d$, using Lemma \ref{A-G-p} we obtain :
$$\dcb
[\widetilde{W^{di}}_j,\widetilde{W^{di}_i}]^{\mathbb{G}-p}

&=&[W^{di}_j,W^{di}_i]^{\mathbb{G}-p}\\
&=&[W^{di}_j,W^{di}_i]^{\mathbb{F}-p}+\Gamma([W^{di}_j,W^{di}_i]-[W^{di}_j,W^{di}_i]^{\mathbb{F}-p})\\

&=&\left((\mathfrak{I}_d+\overline{\mathsf{ci}})\centerdot b\right)_{j,i}\ \ok

\dce
$$

\

\subsection{The main theorem}

In the remainder part of this section we will consider two other equations a little bit stronger than the equations (\ref{structure-condition-da}) and (\ref{structure-condition-di}): For a $\overline{K}\in\underline{\mathsf{L}}^2(\mathbb{P},\mathbb{G},\widetilde{W^{di}}^T)$,
\begin{equation}\label{equation-rho}
\transp\overline{K} \left(\mathfrak{I}_d+\overline{\mathsf{ci}}\right) (\centerdot b)
=(\transp\overline{d}+\transp\overline{\mathsf{f}}) (\centerdot b) \hspace{1cm} \mbox{on $[0,T]$}
\end{equation}
For a $\overline{K}\in\underline{\mathsf{L}}^2(\mathbb{P},\mathbb{G},\widetilde{W^{da}}^T)$,
\begin{equation}\label{equation-kappa}
\transp\overline{K} \left(\mathfrak{I}_d+\overline{\mathsf{ce}} - \Delta a \overline{\mathsf{f}}\transp\overline{\mathsf{f}}\right)(\centerdot a)
=(\transp\overline{d}+\transp\overline{\mathsf{f}}) (\centerdot a) \hspace{1cm} \mbox{on $[0,T]$}
\end{equation}

\brem
We note that these equations (\ref{equation-rho}) and (\ref{equation-kappa}) are not intrinsic versions of the equations (\ref{structure-condition-di}) and (\ref{structure-condition-da}), because they are introduced through the driving process $W$ of a \mrt\ property. We note nevertheless that the equation (\ref{equation-rho}) and (\ref{equation-kappa}) give sufficient conditions to solve the equations (\ref{structure-condition-di}) and (\ref{structure-condition-da}), and we have the freedom to choose suitable driving process $W$ to do the computations. \ok
\erem

Regarding on the equation (\ref{equation-rho}) and (\ref{equation-kappa}), we would attempt to resolve them by assuming that the matrix valued processes $\left(\mathfrak{I}_d+\overline{\mathsf{ci}}\right)$ or $\left(\mathfrak{I}_d+\overline{\mathsf{ce}} - \Delta a \overline{\mathsf{f}}\transp\overline{\mathsf{f}}\right)$ are invertible with bounded inverses. But these conditions are too strong. In fact, we have a more realistic result.

\bassumption\label{1+fin}
\
\ebe
\item[. ]
For any $\mathbb{F}$ predictable stopping time $R$, for any positive random variable $\xi\in\mathcal{F}_R$, we have $\{\mathbb{E}[\xi|\mathcal{G}_{R-}]=0, R\leq T\}=\{\mathbb{E}[\xi|\mathcal{F}_{R-}]=0, R\leq T\}$. 
\item[. ]
There exists a strictly positive $\mathbb{G}$ predictable process $\mathsf{u}$  such that $1+\transp\overline{\varphi}\Delta N\geq \mathsf{u}$ on $[0,T]$.
\dbe
\eassumption

\bethe\label{main}
Suppose Assumptions \ref{assump0}, \ref{assump1} and \ref{assump-mrt}. Suppose Assumption \ref{1+fin} with the process $\mathsf{u}$. Suppose that the increasing processes 
$$
\frac{1}{\mathsf{u}}\centerdot[D^d,D^d] \ \mbox{ and \ } \frac{1}{\mathsf{u}}\transp\overline{\varphi}(\centerdot[N^d,\transp N^d])\overline{\varphi}
$$ 
are $(\mathbb{P},\mathbb{G})$ locally integrable. Then, the equation (\ref{equation-rho}) and the equation (\ref{equation-kappa}), as well as the equations (\ref{structure-condition-da}) and (\ref{structure-condition-di}), have solutions.
\ethe

Applying Lemma \ref{2341} we can now state

\bcor\label{main-consequence}
Under the assumptions in Theorem \ref{main}, if $
\transp\overline{\varphi}(\centerdot[N^c,\transp N^c])\overline{\varphi}
$
defines a finite process, if the increasing process 
$
\frac{1}{\mathsf{u}}\centerdot[D^d,D^d]$ and $\frac{1}{\mathsf{u}}\transp\overline{\varphi}(\centerdot[N^d,\transp N^d])\overline{\varphi}
$
are $(\mathbb{P},\mathbb{G})$ locally integrable, the structure condition (\ref{structure-condition}) holds in the expanded market $(\mathbb{P},\mathbb{G},S)$.
\ecor

The proof of Theorem \ref{main} is developed in the section below.

\

\section{The proof of the main theorem}\label{proof}

Assumptions \ref{assump0}, \ref{assump1} and \ref{assump-mrt} are in force in this section.

\subsection{Conditional multiplicity}\label{multiplicity}

When we try to solve the equations (\ref{equation-rho}) and (\ref{equation-kappa}), because of their jumping nature, we need to compute finely the conditional expectations $\mathbb{E}[\cdot|\mathcal{G}_{R-}]$ for a $\mathbb{F}$ stopping time $R$ (see subsection \ref{exp-at-jump}). These computations will be achieved through the notion of the conditional multiplicity, introduced in \cite[section 3]{BEKSY}, which quantifies the randomness of $\mathcal{F}_{R}$ when $\mathcal{F}_{R-}$ is given.

\bl\label{generating-family}
Let $R$ be a $\mathbb{F}$ stopping time. Consider the random variables in $\mathcal{F}_{R-}$ as constants. If $R$ is predictable, the family of random variables $\Delta_R W_b, 1\leq b\leq d$, generates on $\{R<\infty\}$ (modulo $\mathcal{F}_{R-}$) all integrable random variables $\xi$ in $\mathcal{F}_{R}$ whose conditional expectation $\mathbb{E}[\xi|\mathcal{F}_{R-}]=0$. If $R$ is totally inaccessible, the family of $\Delta_R W_b, 1\leq b\leq d$ generates on $\{R<\infty\}$ (modulo $\mathcal{F}_{R-}$) all integrable random variables $\xi$ in $\mathcal{F}_{R}$.
\el

\textbf{Proof.} For any integrable $\xi\in\mathcal{F}_R$, the process $\xi\ind_{[R,\infty)}-(\xi\ind_{[R,\infty)})^{\mathbb{F}-p}$ is a martingale. By \mrt, there exist $\mathbb{F}$-predictable process $\overline{h}$ such that $\xi\ind_{[T,\infty)}-(\xi\ind_{[T,\infty)})^{\mathbb{F}-p}=\transp\overline{h}\centerdot W$. Therefore, $$
\xi = \sum_{e=1}^{d}(\overline{h}_{R})_e\Delta_R W_e+\Delta_R (\xi\ind_{[R,\infty)})^{\mathbb{F}-p}
$$ 
on $\{R<\infty\}$. If $R$ is predictable and $\mathbb{E}[\xi|\mathcal{F}_{R-}]=0$, $(\xi\ind_{[R,\infty)})^{\mathbb{F}-p}=0$. If $R$ is totally inaccessible, $\Delta_R (\xi\ind_{[R,\infty)})^{\mathbb{F}-p}=0$. The lemma is proved. \ok

\bl\label{partition}
If $R$ is $\mathbb{F}$ predictable, there exists a partition $(A_0,A_1,A_2,\ldots,A_{d})$ (where some $A_i$ may be empty) such that $$
\mathcal{F}_R=\mathcal{F}_{R-}\vee\sigma(A_0,A_1,A_2,\ldots,A_{d})
$$
i.e. the conditional multiplicity of $\mathcal{F}_R$ with respect to $\mathcal{F}_{R-}$ is equal or smaller then $d+1$. If $R$ is $(\mathbb{P},\mathbb{F})$ totally inaccessible, there exists a partition $(B_1,B_2,\ldots,B_{d})$ (where some $B_j$ may be empty) such that $$
\mathcal{F}_R=\mathcal{F}_{R-}\vee\sigma(B_1,B_2,\ldots,B_{d})
$$
i.e. the conditional multiplicity of $\mathcal{F}_R$ with respect to $\mathcal{F}_{R-}$ is equal or smaller then $d$. 
\el

\textbf{Proof.} Consider the case of a predictable $R$. Because of Lemma \ref{generating-family}, we can apply Proposition 12 in \cite{BEKSY} to have a partition $(A'_0,A'_1,A'_2,\ldots,A'_{d})$ of $\{R<\infty\}$ such that $$
\{R<\infty\}\cap\mathcal{F}_R=\{R<\infty\}\cap(\mathcal{F}_{R-}\vee\sigma(A_0,A_1,A_2,\ldots,A_{d})).
$$ 
Since $\{R=\infty\}\cap\mathcal{F}_R=\{R=\infty\}\cap\mathcal{F}_{R-}$, the lemma is verified, if we take $A_i=A'_i$ for $0\leq i<d$ and $A_d=A'_d\cup\{R=\infty\}$.

The case of a totally inaccessible $R$ can be dealt with similarly. \ok

\

\subsection{Martingales with single jump at a $\mathbb{F}$ stopping time}

In this subsection we look at the martingales in the two filtrations, which are the compensated jumps at a $\mathbb{F}$ stopping time.

Consider any $\mathbb{F}$ stopping time $R$. For any real valued $\xi\in\mathbb{L}^1(\mathbb{P},\mathcal{F}_{R})$, let $\ddag_{R}\xi$ denote any coefficient of the $(\mathbb{P},\mathbb{F})$ martingale  $\xi\ind_{[R,\infty)}-(\xi\ind_{[R,\infty)})^{\mathbb{F}-p}$ in its martingale representation with respect to $W$.
(The operation $\ddag$ can be seen as a map whose value is some equivalent class.)

\bl\label{single-jump}
Let $R$ be any $\mathbb{F}$ stopping time either predictable or totally inaccessible. Let $\xi\in\mathbb{L}^1(\mathbb{P},\mathcal{F}_{R})$. If $R$ is predictable, we have
$$
\transp(\ddag_{R}\xi)_{R} \Delta_{R}W=\xi-\mathbb{E}[\xi|\mathcal{F}_{R-}].
$$
If $R$ is totally inaccessible, we have
$$
\transp(\ddag_{R}\xi)_{R} \Delta_{R}W=\xi.
$$
For $\mathbb{F}$ stopping time $S$, we have
$$
\transp(\ddag_{R}\xi)_{S} \Delta_{S}W =0, \mbox{ on $\{S\neq R\}$}.
$$
For any $\mathbb{F}$ predictable precesses $\overline{H}$ in $\underline{\mathsf{L}}(\mathbb{P},\mathbb{F},\xi\ind_{[R,\infty)}-(\xi\ind_{[R,\infty)})^{\mathbb{F}-p})$, $\overline{H}\ddag_{R}\xi$ is in the same equivalence class of $\ddag_{R}(\overline{H_R}\xi)$ in $\underline{\mathsf{L}}(\mathbb{P},\mathbb{F},W)$.
For any $1\leq n<\mathsf{N}^i$ and $1\leq n'<\mathsf{N}^a$, for any $\xi\in\mathbb{L}^1(\mathbb{P},\mathcal{F}_{S_n})$ and $\xi'\in\mathbb{L}^1(\mathbb{P},\mathcal{F}_{T_{n'}})$, 
$$
\dcb
\transp \ddag_{S_n}\xi\centerdot W=\transp \ddag_{S_n}\xi\centerdot W^{di},\hspace{0.5cm} 
\transp \ddag_{T_{n'}}\xi'\centerdot W=\transp \ddag_{T_{n'}}\xi'\centerdot W^{da}\\
\dce
$$
\el

\textbf{Proof.} This lemma is a direct consequence of the properties of stochastic integrals. Let us look at only the first assertions. If $R$ is predictable, $
(\xi\ind_{[R,\infty)})^{\mathbb{F}-p}=\mathbb{E}[\xi|\mathcal{F}_{R-}]\ind_{[R,\infty)}.
$
If $R$ is totally inaccessible, $(\xi\ind_{[R,\infty)})^{\mathbb{F}-p}$ is continuous. Computing the jump at $R$ and at $S$ in the equation$$
\xi\ind_{[R,\infty)}-(\xi\ind_{[R,\infty)})^{\mathbb{F}-p}
=\transp \ddag_{R}\xi\centerdot W
$$
we prove the first assertions. \ok

For $\xi$ a real valued $\mathbb{G}_{R-}$ measurable random variable, we denote by $\vdots \xi\vdots{}^{\mathbb{G}-R}$ the equivalent class of $\mathbb{G}$ predictable processes such that $\vdots \xi\vdots{}^{\mathbb{G}-R}_{R}=\xi$. (We use the same notation to denote also a particular member in the equivalent class) The definition of $\vdots\xi\vdots{}^{\mathbb{G}-R}$ is naturally extended to vector valued $\mathbb{G}_{R-}$ measurable random variables $\xi$. Substituting the filtration, we define also $\vdots \xi\vdots{}^{\mathbb{F}-R}$ for $\mathbb{F}$ stopping times $R$ and $\xi\in\mathcal{F}_{R-}$. The following lemma gives a uniqueness property of the process $\vdots \xi\vdots{}^{\mathbb{G}-R}$.

\bl\label{unique-R}
For any $\mathbb{G}$ stopping time $R$, for any $\mathbb{G}$ predictable process $\overline{H}$, $|\overline{H}|\centerdot(\ind_{[R,\infty)})^{\mathbb{G}-p}=(|\overline{H}_R|\ind_{[R,\infty)})^{\mathbb{G}-p}=0$ if and only if $\overline{H}_R=0$ on $\{R<\infty\}$. The same is true, if the filtration $\mathbb{G}$ is replaced by $\mathbb{F}$.
\el

We consider now the compensated jump in $\mathbb{G}$.

\bl\label{G-jump-R}
For any $\mathbb{F}$ stopping time $R$, for random variables $\xi\in\mathcal{F}_{R}$ and $\xi'\in \mathcal{G}_{R-}$ such that $\xi\in \mathbb{L}^1(\mathbb{P}), \xi'\xi\in\mathbb{L}^1(\mathbb{P})$, we have, on $[0,T]$,$$
\dcb
\xi' \xi\ind_{[R,\infty)}-(\xi'\xi\ind_{[R,\infty)})^{\mathbb{G}-p}
=\vdots \xi'\vdots{}^{\mathbb{G}-R}\transp (\ddag_{R} \xi)\centerdot\widetilde{W}
\dce
$$
\el

\textbf{Proof.} We need only to note that, by Lemma \ref{A-G-p}, on $[0,T]$,$$
\dcb
\xi\ind_{[R,\infty)}-(\xi\ind_{[R,\infty)})^{\mathbb{G}-p}
=\xi\ind_{[R,\infty)}-(\xi\ind_{[R,\infty)})^{\mathbb{F}-p}-\Gamma(\xi\ind_{[R,\infty)}-(\xi\ind_{[R,\infty)})^{\mathbb{F}-p})

&=&\transp (\ddag_{R} \xi)\centerdot\widetilde{W}\ \ok
\dce
$$

\

\

\subsection{Equations at the stopping times $S_n$ or $T_{n'}$}

Before solving the equations (\ref{equation-rho}) and (\ref{equation-kappa}) we note that these equations impose a strict relationship between the jumps in the two filtrations $\mathbb{F}$ and $\mathbb{G}$.

Let $(S_n)_{1\leq n<\mathsf{N}^i}$ ($\mathsf{N}^i\leq \infty$) (resp. $(T_n)_{1\leq n<\mathsf{N}^a}$) be a sequence of $(\mathbb{P},\mathbb{F})$ totally inaccessible (resp. strictly positive $(\mathbb{P},\mathbb{F})$ predictable) stopping times such that $[S_n]\cap [S_{n'}]=\emptyset$ for $n\neq n'$ and $\{s\geq 0:\Delta_sW^{di}\neq 0\}\subset\cup_{n\geq 1}[S_n]$ (resp. $[T_n]\cap [T_{n'}]=\emptyset$ for $n\neq n'$ and $\{s\geq 0:\Delta_sW^{da}\neq 0\}\subset\cup_{n\geq 1}[T_n]$).

Consider firstly the case of the predictable stopping times $T_n$.

\bl\label{at-Tn}
On $\{T_n\leq T, T_n<\infty\}$, we have $$
\mathbb{E}[\Delta_{T_n}\widetilde{W}\transp\Delta_{T_n}\widetilde{W}|\mathcal{G}_{T_n-}]=\left(\mathfrak{I}_d+\overline{\mathsf{ce}} - \Delta a \overline{\mathsf{f}}\transp\overline{\mathsf{f}}\right)_{T_n}\Delta_{T_n}a
$$ 
In particular, $\left(\mathfrak{I}_d+\overline{\mathsf{ce}} - \Delta a \overline{\mathsf{f}}\transp\overline{\mathsf{f}}\right)\Delta a$ is symmetric.
If $\overline{K}$ satisfies the equation (\ref{equation-kappa}), then, for every $1\leq n<\mathsf{N}^a$, $\overline{K}_{T_n}$ satisfies the equation $$
\transp\overline{K}_{T_n} \left(\mathfrak{I}_d+\overline{\mathsf{ce}} - \Delta a \overline{\mathsf{f}}\transp\overline{\mathsf{f}}\right)_{T_n}\Delta_{T_n}a
=(\transp\overline{d}+\transp\overline{\mathsf{f}})_{T_n} \Delta_{T_n}a
$$
or equivalently
$$
\transp\overline{K}_{T_n} \mathbb{E}[\Delta_{T_n}\widetilde{W}\transp \Delta_{T_n}\widetilde{W}|\mathcal{G}_{T_n-}]
=(\transp\overline{d}+\transp\overline{\mathsf{f}})_{T_n}\mathbb{E}[\Delta_{T_n}W\transp \Delta_{T_n}W|\mathcal{F}_{T_n-}] 
$$
on $\{T_n\leq T, T_n<\infty\}$.
\el

\textbf{Proof.} The first assertion is a direct consequence of the formula (\ref{jcef}). To prove the other assertion, we need only to compute the jump at $T_n$ in the equation (\ref{equation-kappa}). \ok

Consider now the case of the totally inaccessible stopping times $S_n$. We introduce the notion $g_{S_n}$ to design the $d\times d$ matrix valued random variable whose coefficients are defined by$$
(g_{S_n})_{j,i}=\mathbb{E}[\Delta_{S_n} W_j\Delta_{S_n} W_i|\mathcal{F}_{S_n-}],\ 1\leq i,j\leq d.
$$

\bl\label{individual-equation}
If $\overline{K}$ satisfies the equation (\ref{equation-rho}), it satisfies also, at every $S_n,1\leq n<\mathsf{N}^i$, $$
\transp\overline{K}(\mathfrak{I}_d+\overline{\mathsf{ci}})\centerdot (g_{S_n}\ind_{[S_n,\infty)})^{\mathbb{F}-p}
=(\transp\overline{d}+\transp\overline{\mathsf{f}})\centerdot (g_{S_n}\ind_{[S_n,\infty)})^{\mathbb{F}-p}\ \mbox{ on $[0,T]$}
$$
\el

\textbf{Proof.} For $1\leq h\leq d$, integrating $(\ddag_{S_n} \Delta_{S_n}W_h)$ into the two side of the equation (\ref{equation-rho}) we obtain :$$
\transp\overline{K} \left(\mathfrak{I}_d+\overline{\mathsf{ci}} \right)(\centerdot b) (\ddag_{S_n} \Delta_{S_n}W_h)
=(\transp\overline{d}+\transp\overline{\mathsf{f}}) (\centerdot b) (\ddag_{S_n} \Delta_{S_n}W_h)\
$$
on $[0,T]$. The lemma follows, because$$
\dcb
(\centerdot b)(\ddag_{S_n} \Delta_{S_n}W_h)

&=&((g_{S_n})_{\cdot,h}\ind_{[S_n,\infty)})^{\mathbb{F}-p}\ \ok
\dce
$$

\bl\label{varphi-r}
Let $1\leq n<\mathsf{N}^i$. Let $\mathtt{r}_n$ denote a vector valued process in $\vdots\mathbb{E}[\Delta_{S_n}N|\mathcal{F}_{S_n-}]\vdots{}^{\mathbb{F}-S_n}$. Then, $$
(\ind_{[S_n,\infty)})^{\mathbb{G}-p}
=\left(1+\transp\overline{\varphi}\mathtt{r}_n\right)\centerdot(\ind_{[S_n,\infty)})^{\mathbb{F}-p}
$$ 
on $[0,T]$. Consequently, $(\ind_{[S_n,\infty)})^{\mathbb{G}-p}$ is continuous (i.e. $S_n$ is $(\mathbb{P},\mathbb{G})$ totally inaccessible), and $1+\overline{\varphi}\mathtt{r}_n\geq 0$ on $[0,T]$ almost surely under the random measure $\mathsf{d}(\ind_{[S_n,\infty)})^{\mathbb{F}-p}$. In particular, for any $\xi\in\mathcal{F}_{S_n-}$, $\vdots\xi\vdots{}^{\mathbb{F}-S_n}\subset\vdots\xi\vdots{}^{\mathbb{G}-S_n}$ on $[0,T]$.
\el

\textbf{Proof.} The first assertion is the consequence of the following computations from Lemma \ref{A-G-p} on $[0,T]$ :
$$
\dcb
&&
(\ind_{[S_n,\infty)})^{\mathbb{G}-p}\\

&=&
(\ind_{[S_n,\infty)})^{\mathbb{F}-p}+\transp\overline{\varphi}_{S_n}(\Delta_{S_n}N\ind_{[S_n,\infty)})^{\mathbb{F}-p}\\

&=&
\left(1+\transp\overline{\varphi}\mathtt{r}_n\right)\centerdot(\ind_{[S_n,\infty)})^{\mathbb{F}-p}\ \ok
\dce
$$
Since $(\ind_{[S_n,\infty)})^{\mathbb{G}-p}$ and $(\ind_{[S_n,\infty)})^{\mathbb{F}-p}$ are increasing, we get the second assertion. The other assertions of the lemma are obvious.\ok

Combining the preceding lemmas, we obtain the equations at the stopping times $S_n$:

\bcor\label{at-Sn}
Let $1\leq n<\mathsf{N}^i$. If $\overline{K}$ satisfies the equation (\ref{equation-rho}), then $\overline{K}_{S_n}$ satisfies the equation $$
\transp\overline{K}_{S_n}(\mathfrak{I}_d+\overline{\mathsf{ci}})_{S_n}g_{S_n}
= (\transp\overline{d}+\transp\overline{\mathsf{f}})_{S_n}g_{S_n}\ \mbox{ on $\{S_n\leq T, S_n<\infty\}$}
$$
\ecor

\textbf{Proof.} According to Lemma \ref{individual-equation}, we have on $[0,T]$$$
(1+\transp\overline{\varphi}\mathtt{r}_n)\transp\overline{K}(\mathfrak{I}_d+\overline{\mathsf{ci}}) \centerdot(g_{S_n}\ind_{[S_n,\infty)})^{\mathbb{F}-p}
=
(1+\transp\overline{\varphi}\mathtt{r}_n)(\transp\overline{d}+\transp\overline{\mathsf{f}})\centerdot (g_{S_n}\ind_{[S_n,\infty)})^{\mathbb{F}-p}
$$
Using Lemma \ref{varphi-r}, we compute on $[0,T]$, $$
\dcb
(1+\transp\overline{\varphi}\mathtt{r}_n)\transp\overline{K}(\mathfrak{I}_d+\overline{\mathsf{ci}}) \centerdot(g_{S_n}\ind_{[S_n,\infty)})^{\mathbb{F}-p}

&=&(\transp\overline{K}_{S_n}(\mathfrak{I}_d+\overline{\mathsf{ci}})_{S_n} g_{S_n}\ind_{[S_n,\infty)})^{\mathbb{G}-p}\\
\dce
$$
and
$$
\dcb
(1+\transp\overline{\varphi}\mathtt{r}_n)(\transp\overline{d}+\transp\overline{\mathsf{f}})\centerdot (g_{S_n}\ind_{[S_n,\infty)})^{\mathbb{F}-p}

&=&((\transp\overline{d}+\transp\overline{\mathsf{f}})_{S_n} g_{S_n}\ind_{[S_n,\infty)})^{\mathbb{G}-p}
\dce
$$
Comparing these two identities, applying Lemma \ref{unique-R}, we prove the lemma. \ok

We can express the equation at $S_n$ in term of conditional expectations. We begin with the following lemma that can be proved with a direct computation.

\bl\label{cig}
We have$$
\sum_{k=1}^n(\overline{\varphi}_{S_n})_k\mathbb{E}[\Delta_{S_n}N_k\Delta_{S_n}W\transp \Delta_{S_n}W|\mathcal{F}_{S_n-}]
=\overline{\mathsf{ci}}_{S_n}g_{S_n}
$$
\el

\bl\label{conditional-expectation-i}
Let $1\leq n<\mathsf{N}^i$. Let $\mathtt{r}_n$ be the vector valued process introduced in Lemma \ref{varphi-r}. We have 
$$
(1+\transp\overline{\varphi}_{S_n}\mathtt{R}_n)\mathbb{E}[\Delta_{S_n}W\transp \Delta_{S_n}W|\mathcal{G}_{S_n-}]
=(\mathfrak{I}_d+\overline{\mathsf{ci}})_{S_n} g_{S_n}
$$
on $\{S_n\leq T, S_n<\infty\}$, with $\mathtt{R}_{n}=(\mathtt{r}_n)_{S_n}=\mathbb{E}[\Delta_{S_n}N|\mathcal{F}_{S_n-}]$. In particular, $(\mathfrak{I}_d+\overline{\mathsf{ci}})_{S_n} g_{S_n}$ is a symmetric matrix there.
\el

\textbf{Proof.} Applying Lemmas \ref{A-G-p}, \ref{cig} and \ref{varphi-r}, we compute on $[0,T]$
$$
\dcb
&&
((1+\transp\overline{\varphi}_{S_n}\mathtt{R}_n)\Delta_{S_n}W\transp \Delta_{S_n}W\ind_{[S_n,\infty)})^{\mathbb{G}-p}\\

&=&
(1+\transp\overline{\varphi}\mathtt{r}_n)\centerdot\left((\Delta_{S_n}W\transp \Delta_{S_n}W\ind_{[S_n,\infty)})^{\mathbb{F}-p}+\sum_{k=1}^n\overline{\varphi}_k\centerdot(\Delta_{S_n}N_k\Delta_{S_n}W\transp \Delta_{S_n}W\ind_{[S_n,\infty)})^{\mathbb{F}-p}\right)\\

&=&
((\mathfrak{I}_d+\overline{\mathsf{ci}})_{S_n}g_{S_n}\ind_{[S_n,\infty)})^{\mathbb{G}-p}\\
\dce
$$
This computation together with Lemma \ref{unique-R} prove the result.\ok

\bcor\label{KWG}
Let $1\leq n<\mathsf{N}^i$. The equation in Corollary \ref{at-Sn} have an expression in term of conditional expectations :
$$
(1+\transp\overline{\varphi}_{S_n}\mathtt{R}_n)\transp\overline{K}\mathbb{E}[\Delta_{S_n}W\transp \Delta_{S_n}W|\mathcal{G}_{S_n-}]
=(\transp\overline{d}+\transp\overline{\mathsf{f}})_{S_n}\mathbb{E}[\Delta_{S_n}W\transp \Delta_{S_n}W|\mathcal{F}_{S_n-}]\ 
$$
on $\{S_n\leq T, S_n<\infty\}$
\ecor

\

\subsection{Constructions of solutions}

Now a natural idea to solve the equations (\ref{equation-rho}) and (\ref{equation-kappa}) is to solve firstly the equations in Lemma \ref{at-Tn} and in Corollary \ref{at-Sn} at $T_{n'}$ and at $S_n$, and then to integrate the individual solutions with each other to form a global solution.

We need some algebraic properties.

\bl\label{symmetric-matrix}
Let $G$ and $J$ be two $d\times d$ matrix such that $G$ and $GJ$ are symmetric and are positive semidefinite. We consider $G,J$ as linear operators on $\mathbb{R}^d$. We identify the matrix $G$ with the map $v\rightarrow Gv, v\in\mathbb{R}^d$. Same for $J$. Let $\mathbb{V}$ be the image space of $G$. Recall that $(\!(v|v)\!), v\in \mathbb{R}^d$, denotes the inner product in $\mathbb{R}^d$. Then,
\ebe
\item[. ]
If $G$ is not trivial, $G$ as a linear map is invertible on $\mathbb{V}$
\item[. ]
$\transp J$ is symmetric with respect to the quadratic form $(\!(v|Gv)\!), v\in\mathbb{R}^d$.
\item[. ]
If $G$ is not trivial, if there exists a $\epsilon>0$ such that $(\!(v|GJ v)\!)\geq \epsilon (\!(v|Gv)\!), \forall v\in\mathbb{R}^d$, then, $J$ is an invertible operator on $\mathbb{V}$. Let $J^*$ denote its inverse on $\mathbb{V}$, we have $$
(\!(J^*v|GJ^* v)\!)^{1/2}\leq \frac{1}{\epsilon} (\!(v|Gv)\!)^{1/2}, \forall v\in\mathbb{R}^d
$$
When $G$ and $J$ are measurable functions on some space, $J^*$ is also a measurable function.
\item[. ]
Let $\mathtt{p}_G$ be the orthogonal projection onto $\mathbb{V}$ (orthogonality with respect to the inner product $(\!(\cdot|\cdot)\!)$ in $\mathbb{R}^d$). We have $G\mathtt{p}_Gv=Gv, v\in\mathbb{R}^d$. If $G$ is a measurable function on some space, $\mathtt{p}_G$ is also a measurable function.
\dbe
\el

The following lemma constructs a solution to the equation (\ref{equation-kappa}).

\bl\label{da-xi-lemma}
Suppose that, for any $1\leq n<\mathsf{N}^a$, there exists a strictly positive random variable $\varkappa_n$ such that 
\begin{equation}\label{da-u-inequality}
(\!(v|(\Delta_{T_n}a)\transp(\mathfrak{I}_d+\overline{\mathsf{ce}}-\Delta a\overline{\mathsf{f}}\transp\overline{\mathsf{f}})_{T_n} v)\!)\
\geq \varkappa_n\ (\!(v|(\Delta_{T_n}a) v)\!),\ \forall v=(v_h)_{1\leq h\leq d}\in\mathbb{R}^d
\end{equation}
on $\{T_n\leq T, T_n<\infty\}$. Let $$
\xi_n=\ind_{\{T_n\leq T, T_n<\infty\}}\ind_{\{\Delta_{T_n}a\neq 0\}}\left(\transp(\mathfrak{I}_d+\overline{\mathsf{ce}}-\Delta a\overline{\mathsf{f}}\transp\overline{\mathsf{f}})_{T_n}\right)^*\mathtt{p}_{\Delta_{T_n}a}(\overline{d}+\transp\overline{\mathsf{f}})_{T_n}
$$ 
which is a $d$-dimensional vector valued $\mathbb{G}_{T_n-}$ measurable random variable. Suppose that the increasing process $
\sum_{1\leq n<\mathsf{N}^a}(\transp\xi_{n}\Delta_{T_n}\widetilde{W})^2\ind_{[T_n,\infty)}
$
is $(\mathbb{P},\mathbb{G})$ locally integrable. Then, the series $\overline{K}=\sum_{1\leq n<\mathsf{N}^a}\xi_n\ind_{[T_n]}$ converges in the space $\underline{\mathsf{L}}^2(\mathbb{P},\mathbb{G},\widetilde{W^{da}}^T)$ and the equation (\ref{equation-kappa}) is solved by $\overline{K}$.
\el

\textbf{Proof.} The random variables $\xi_n$ are well defined according to Lemma \ref{symmetric-matrix}. To see the convergence of the series $\overline{K}$, let $(R_m)$ be an increasing sequence of $\mathbb{G}$ stopping times such that $\sup_nR_n=T$ and, for every $n\geq 1$, the random variable $
\sum_{1\leq n<\mathsf{N}^a}(\transp\xi_{n}\Delta_{T_n}\widetilde{W})^2\ind_{\{T_n\leq R_m\}}
$
is $\mathbb{P}$ integrable. We have the following inequalities $$
\dcb
&&\transp(\sum_{i\leq n\leq j}\xi_n\ind_{[T_n]})(\centerdot [\widetilde{W},\transp \widetilde{W}])(\sum_{i\leq n\leq j}\xi_n\ind_{[T_n]})_{R_m}\\
&=&[(\transp(\sum_{i\leq n\leq j}\xi_n\ind_{[T_n]})\centerdot \widetilde{W},(\transp(\sum_{i\leq n\leq j}\xi_n\ind_{[T_n]})\centerdot \widetilde{W}]_{R_m}\\

&=&\sum_{i\leq n\leq j}(\transp\xi_{n}\Delta_{T_n}\widetilde{W})^2\ind_{\{T_n\leq R_m\}}\\
\dce
$$
The last term converges in $\mathbb{L}^1(\mathbb{P})$ to zero when $i,j\uparrow\infty$, which proves the convergence of $\overline{K}$. 

To check that $\overline{K}$ solves the equation (\ref{equation-kappa}), we note that
$\xi_n$ satisfies the equations : 
$$
\transp\xi_{n}(\mathfrak{I}_d+\overline{\mathsf{ce}}-\Delta a\overline{\mathsf{f}}\transp\overline{\mathsf{f}})_{T_n}\Delta_{T_n}a
=
(\transp\overline{d}+\transp\overline{\mathsf{f}})_{T_n}\Delta_{T_n}a
$$
on $\{T_n\leq T, T_n<\infty\}$, and consequently, on $[0,T]$$$
\dcb
\transp\overline{K}\centerdot(\mathfrak{I}_d+\overline{\mathsf{ce}}-\Delta a\overline{\mathsf{f}}\transp\overline{\mathsf{f}})\centerdot a

=
\sum_{1\leq n<\mathsf{N}^a}\transp\xi_n(\mathfrak{I}_d+\overline{\mathsf{ce}}-\Delta a\overline{\mathsf{f}}\transp\overline{\mathsf{f}})_{T_n}\Delta_{T_n} a

&=&
(\transp\overline{d}+\transp\overline{\mathsf{f}})\centerdot a \ \ok
\dce
$$

We construct now a solution to the equation (\ref{equation-rho}). We recall the decompositions of the jumping parts of $W$ (see \cite{HWY}). 

\bl\label{W-jumps}
The jumping part of the driving process $W$ is given by $$
\dcb
W^{di}=\sum_{1\leq n<\mathsf{N}^i}(\Delta_{S_n}W^{di}\ind_{[S_n,\infty)}-(\Delta_{S_n}W^{di}\ind_{[S_n,\infty)})^{\mathbb{F}-p})\\
W^{da}=\sum_{1\leq n<\mathsf{N}^a}\Delta_{T_n}W^{da}\ind_{[T_n,\infty)}
\dce
$$
where the series converge in $\mathcal{H}^2_{loc}(\mathbb{P},\mathbb{F})$. Consequently,$$
\dcb
b=[W^{di},\transp W^{di}]^{\mathbb{F}-p}
=\sum_{1\leq n<\mathsf{N}^i}\left(g_{S_n}\ind_{[S_n,\infty)}\right)^{\mathbb{F}-p}\\ 
a=[W^{da},W^{da}]^{\mathbb{F}-p}
=\sum_{1\leq n<\mathsf{N}^a}\Delta_{T_n}a\ind_{[T_n,\infty)}
\dce
$$
\el

\bl\label{di-xi-lemma}
Suppose that, for $1\leq n<\mathsf{N}^i$, there exists a strictly positive random variable $\varrho_n$ such that 
\begin{equation}\label{di-u-inequality}
(v|g_{S_n}\transp(\mathfrak{I}_d+\overline{\mathsf{ci}})_{S_n} v)\
\geq \varrho_n\ (v|g_{S_n} v),\ \forall v=(v_h)_{1\leq h\leq d}\in\mathbb{R}^d
\end{equation}
on $\{S_n\leq T, S_n<\infty\}$. Let $$
\xi_n=\ind_{\{S_n\leq T, S_n<\infty\}}\ind_{\{g_{S_n}\neq 0\}}\left(\transp(\mathfrak{I}_d+\overline{\mathsf{ci}})_{S_n}\right)^*\mathtt{p}_{g_{S_n}}(\overline{d}+\overline{\mathsf{f}})_{S_n}
$$ 
which is a $d$-dimensional vector valued $\mathbb{G}_{S_n-}$ measurable random variable. Suppose that the increasing process $
\sum_{1\leq n<\mathsf{N}^i}(\transp\xi_{n}\Delta_{S_n}W)^2\ind_{[S_n,\infty)}
$
is $(\mathbb{P},\mathbb{G})$ locally integrable. Define $\ddag_{S_n} \Delta_{S_n}W$ to be the matrix valued process whose $h$-line vector is given by $\transp\ddag_{S_n} \Delta_{S_n}W_h$. Define $\overline{K}_n$ to be the vector valued process  
$$
\overline{K}_{n}=\transp\vdots\xi_{n}\vdots{}^{\mathbb{G}-S_n}\ddag_{S_n} \Delta_{S_n}W,
$$ 
for $1\leq h\leq d$. Then, the series $\overline{K}=\sum_{1\leq n<\mathsf{N}^i}\overline{K}_n$ converges in the space $\underline{\mathsf{L}}^2(\mathbb{P},\mathbb{F},\widetilde{W^{di}}^T)$ and the equation (\ref{equation-rho}) is solved by $\overline{K}$.
\el

\textbf{Proof.} The random variables $\xi_n$ are well defined according to Lemma \ref{symmetric-matrix}. From Lemma \ref{G-jump-R} and \ref{single-jump}, we can write on $[0,T]$ $$
\overline{K}_n\centerdot\widetilde{W^{di}}=\transp\vdots\xi_{n}\vdots{}^{\mathbb{G}-S_n} \ddag_{S_n} \Delta_{S_n}W\centerdot\widetilde{W}
=\transp\xi_{n} \Delta_{S_n}W\ind_{[S_n,\infty)}-(\transp\xi_{n}\Delta_{S_n}W\ind_{[S_n,\infty)})^{\mathbb{G}-p}
$$ 
It follows that, for $1\leq m<m'$:$$
\dcb
[\transp (\sum_{n=m}^{m'}\overline{K}_{n})\centerdot\widetilde{W^{di}},\transp (\sum_{n=m}^{m'}\overline{K}_{n})\centerdot\widetilde{W^{di}}]
&=&\sum_{n=m}^{m'}(\transp\xi_{n} \Delta_{S_n}W)^2 \ind_{[S_n,\infty)}.
\dce
$$
The assumption of the lemma implies then the convergence of $\overline{K}$. It is to note then $$
\transp\overline{K}\centerdot \widetilde{W^{di}}=\sum_{1\leq n<\mathsf{N}^i}\transp\overline{K}_n\centerdot \widetilde{W^{di}}
=\sum_{1\leq n<\mathsf{N}^i}(\xi_{n} \Delta_{S_n}W\ind_{[S_n,\infty)}-(\xi_{n}\Delta_{S_n}W\ind_{[S_n,\infty)})^{\mathbb{G}-p})
$$
on $[0,T]$, where the series converges in $\mathcal{H}^2_{loc}(\mathbb{P},\mathbb{G})$. 

Now to see that $\overline{K}$ solves the equation (\ref{equation-rho}), we note that $
\transp\xi_{n}(\mathfrak{I}_d+\overline{\mathsf{ci}})_{S_n}g_{S_n}=(\transp\overline{d}+\transp\overline{\mathsf{f}})_{S_n}g_{S_n},
$
or, according to Lemma \ref{conditional-expectation-i}, 
$$
(1+\transp\overline{\varphi}_{S_n}\mathtt{R}_n)\transp\xi_n\mathbb{E}[\Delta_{S_n}W\transp \Delta_{S_n}W|\mathcal{G}_{S_n-}]
=\transp\xi_n(\mathfrak{I}_d+\overline{\mathsf{ci}})_{S_n} g_{S_n}
=(\transp\overline{d}+\transp\overline{\mathsf{f}})_{S_n}g_{S_n}
$$
Applying then Lemma \ref{sc-di}, Lemma \ref{varphi-r} and Lemma \ref{F-Sn} (by anticipation) we compute on $[0,T]$$$
\dcb
&&\transp\overline{K}(\mathfrak{I}_d+\overline{\mathsf{ci}})(\centerdot b)\\

&=&\transp\overline{K}\centerdot[\widetilde{W^{di}},\transp \widetilde{W^{di}}]^{\mathbb{G}-p}\\

&=&
\sum_{1\leq n<\mathsf{N}^i}(\transp\xi_{n}\Delta_{S_n}W\transp \Delta_{S_n}W\ind_{[S_n,\infty)})^{\mathbb{G}-p}\\

&=&
\sum_{1\leq n<\mathsf{N}^i}\frac{1}{(1+\transp\overline{\varphi}\mathtt{r}_n)}(\transp\overline{d}+\transp\overline{\mathsf{f}})\centerdot(g_{S_n} \ind_{[S_n,\infty)})^{\mathbb{G}-p}\\

&=&
(\transp\overline{d}+\transp\overline{\mathsf{f}})\centerdot b\ \ok
\dce
$$

\subsection{Conditional expectations at the jumping times}\label{exp-at-jump}

In this subsection, applying the notion of conditional multiplicity (see subsection \ref{multiplicity}), we compute finely the conditional expectations at $S_n$ or at $T_{n'}$. The results of this subsection will be crucial to make use of Lemma \ref{da-xi-lemma} and \ref{di-xi-lemma}.

\bl\label{F-Tn}
For a fixed $1\leq n<\mathsf{N}^a$, let $(A_0,A_1,A_2,\ldots,A_{d})$ be a partition which satisfies the relation $
\mathcal{F}_{T_n}=\mathcal{F}_{T_n-}\vee\sigma(A_0,A_1,A_2,\ldots,A_{d}).
$
Denote $p_k=\mathbb{P}[A_k|\mathcal{F}_{T_n-}], 0\leq k\leq d$.  
\ebe
\item[. ]
For any finite random variable $\xi\in\mathcal{F}_{T_n}$, the conditional expectation $\mathbb{E}[\xi|\mathcal{F}_{T_n-}]$ is well-defined. Let $$
\mathfrak{a}(\xi)_k=\ind_{\{p_k>0\}}\frac{1}{p_k}\mathbb{E}[\ind_{A_k}\xi|\mathcal{F}_{T_n-}], \ 0\leq k\leq d.
$$
We have $\xi=\sum_{h=0}^d\mathfrak{a}(\xi)_h\ind_{A_h}$. We can consider $\xi$ as a random function define on the measurable space $\{0,1,2,\ldots,d\}$ (the value of $\xi$ at point $k$ being $\mathfrak{a}(\xi)_k$). 
\item[. ]
Denote $n$-dimensional vector valued random variable $n_k=\mathfrak{a}(\Delta_{T_n}N)_k, 0\leq k\leq d$. We have $
(1+\transp\overline{\varphi}_{T_n}n_k)p_k
=
\mathbb{E}[\ind_{A_k}|\mathcal{G}_{T_n-}]\geq 0
$
for $0\leq k\leq d$, on $\{T_n\leq T, T_n<\infty\}$.
This means that, if $\chi_{T_n}$ denotes the random variable defined by $\chi_{T_n}=\sum_{k=0}^dk\ind_{A_k}\in\{0,1,2,\ldots,d\}$,
the conditional law of $\chi_{T_n}$ on $\{T_n\leq T, T_n<\infty\}$ given $\mathcal{G}_{T_n-}$ is absolutely continuous with respect to the conditional law of $\chi_{T_n}$ given $\mathcal{F}_{T_n-}$ with a (random) density function $\sum_{k=0}^d(1+\transp\overline{\varphi}_{T_n}n_k)\ind_{\{k\}}$ 
  
\item[. ]
We have the identity :$$
\mathbb{E}[(\!(v|\Delta_{T_n}\widetilde{W})\!)^2|\mathcal{G}_{T_n-}]
=
\sum_{h=0}^d(1+(\transp\overline{\varphi}_{T_n}n_h))p_h\left((\!(v|w_h)\!)-\sum_{k=0}^d(1+(\transp\overline{\varphi}_{T_n}n_k))(\!(v|w_k)\!)p_k\right)^2
$$
for $v\in\mathbb{R}^d$, $w_k=\mathfrak{a}(\Delta_{T_n}W)_k, 0\leq k\leq d$, on $\{T_n\leq T, T_n<\infty\}$.
\dbe
\el

\textbf{Proof.}
The first assertion of the lemma is the consequence of the relation $
\mathcal{F}_{T_n}=\mathcal{F}_{T_n-}\vee\sigma(A_0,A_1,A_2,\ldots,A_{d}).
$
The second assertion follows from a direct computation using Lemma \ref{A-G-p}.

The probability density property of $(1+\transp\overline{\varphi}_{T_n}n_k)\ind_{\{k\}}$ is the consequence of the second assertion and of the vanishing identity $\sum_{k=0}^d \transp\overline{\varphi}_{T_n}n_kp_k=0$, because $N$ is a $(\mathbb{P},\mathbb{F})$ local martingale. As for the last assertion, we recall the formula in Lemma \ref{at-Tn} and, for $1\leq j,i\leq d$,$$
\dcb
\sum_{e=1}^d(\overline{\mathsf{ce}}_{T_n})_{j,e}\mathbb{E}[\Delta_{T_n}W_e\Delta_{T_n}W_{i}|\mathcal{F}_{T_n-}]

&=&\transp\overline{\varphi}_{T_n}\mathbb{E}[\Delta_{T_n}N\Delta_{T_n}W_{j}\Delta_{T_n}W_{i}|\mathcal{F}_{T_n-}]\\

\dce
$$
Now we compute on $\{T_n\leq T, T_n<\infty\}$ for $v\in\mathbb{R}^d$$$
\dcb
&&\mathbb{E}[(\transp v\Delta_{T_n}\widetilde{W})^2|\mathcal{G}_{T_n-}]\\

&=&(\!(v|(\mathfrak{I}_d+\overline{\mathsf{ce}}-\Delta a\overline{\mathsf{f}}\transp\overline{\mathsf{f}})_{T_n} (\Delta_{T_n}a)v)\!)\\

&=&\mathbb{E}[(\!(v|\Delta_{T_n}W)\!)^2|\mathcal{F}_{T_n-}]
+\transp\overline{\varphi}_{T_n}\mathbb{E}[\Delta_{T_n}N(\!(v|\Delta_{T_n}W)\!)^2|\mathcal{F}_{T_n-}]
-\left(\transp\overline{\varphi}_{T_n}\mathbb{E}[\Delta_{T_n}N(\!(v|\Delta_{T_n}W)\!)|\mathcal{F}_{T_n-}]\right)^2\\

&=&\sum_{h=0}^d(\!(v|w_h)\!)^2p_h
+\transp\overline{\varphi}_{T_n}\sum_{h=0}^dn_h(\!(v|w_h)\!)^2p_h
-\left(\transp\overline{\varphi}_{T_n}\sum_{h=0}^dn_h(\!(v|w_h)\!)p_h\right)^2\\

&=&\sum_{h=0}^d(1+(\transp\overline{\varphi}_{T_n}n_h))p_h\left((\!((v|w_h)\!)-(\sum_{k=0}^d(1+(\transp\overline{\varphi}_{T_n}n_k))(\!(v|w_k)\!)p_k)\right)^2\ \ok
\dce
$$

\bcor\label{w-daf}
Let $1\leq n<\mathsf{N}^a$. We use the same notations as in the preceding lemma. If $\overline{K}$ satisfies the equation (\ref{equation-kappa}), then, for $0\leq k\leq d$, $$
(\transp\overline{K}_{T_n}(w_k-\Delta_{T_n}a \overline{\mathsf{f}})-1) (1+\transp\overline{\varphi}_{T_n}n_k)p_k
=(d_k-1)p_k \leq 0
$$
on $\{T_n\leq T, T_n<\infty\}$, and $\transp\overline{K}_{T_n}(w_k-\Delta_{T_n}a \overline{\mathsf{f}})<1$, $(1+\transp\overline{\varphi}_{T_n}n_k)>0$ if and only if $p_k>0$.
\ecor

\textbf{Proof.} According to Lemma \ref{at-Tn}, we can write on $\{T_n\leq T, T_n<\infty\}$:
$$
\transp\overline{K}_{T_n} \mathbb{E}[\Delta_{T_n}\widetilde{W}\transp \Delta_{T_n}\widetilde{W}|\mathcal{G}_{T_n-}](\ddag_{T_n}\ind_{A_k})_{T_n}
=(\transp\overline{d}+\transp\overline{\mathsf{f}})_{T_n}\mathbb{E}[\Delta_{T_n}W\transp \Delta_{T_n}W|\mathcal{F}_{T_n-}] (\ddag_{T_n}\ind_{A_k})_{T_n}
$$
on $\{T_n\leq T, T_n<\infty\}$. Applying Lemma \ref{single-jump} and Lemma \ref{G-jump-R}, we write 
$$
\transp\overline{K}_{T_n} \mathbb{E}[\Delta_{T_n}\widetilde{W} \ind_{A_k}|\mathcal{G}_{T_n-}]
=(\transp\overline{d}+\transp\overline{\mathsf{f}})_{T_n}\mathbb{E}[\Delta_{T_n}W \ind_{A_k}|\mathcal{F}_{T_n-}] 
$$
But$$
\dcb
\Delta_{T_n}\widetilde{W} \ind_{A_k}

&=&(w_k-\Delta_{T_n}a \overline{\mathsf{f}})\ind_{A_k}
\dce
$$
and
$$
(\transp\overline{d}+\transp\overline{\mathsf{f}})_{T_n}\Delta_{T_n}W \ind_{A_k}
=
\Delta_{T_n}D \ind_{A_k}+\transp\overline{\varphi}_{T_n}\Delta_{T_n}N \ind_{A_k}
=
(d_k+\transp\overline{\varphi}_{T_n}n_k) \ind_{A_k}
$$
where $d_k=\mathfrak{a}(\Delta_{T_n}D)_k$. It follows that
$$
\transp\overline{K}_{T_n}(w_k-\Delta_{T_n}a \overline{\mathsf{f}}) \mathbb{E}[\ind_{A_k}|\mathcal{G}_{T_n-}]
=(d_k+\transp\overline{\varphi}_{T_n}n_k)\mathbb{E}[ \ind_{A_k}|\mathcal{F}_{T_n-}] 
$$
or equivalently
$$
(\transp\overline{K}_{T_n}(w_k-\Delta_{T_n}a \overline{\mathsf{f}})-1) (1+\transp\overline{\varphi}_{T_n}n_k))p_k
=(d_k-1)p_k 
$$
Note that, by Assumption \ref{assump0}, $\Delta_{T_n}D-1<0$ which yields $(d_k-1)p_k\leq 0$. Note equally that $$
0=\ind_{\{d_h-1=0\}}(d_h-1)\ind_{A_h}=\ind_{\{d_h-1=0\}}(\Delta_{T_n}N-1)\ind_{A_h}.
$$
This means that $\ind_{\{d_h-1=0\}}\ind_{A_h}=0$. Taking conditional expectation with respect to $\mathcal{F}_{T_n-}$ we have also $\ind_{\{d_h-1=0\}}p_h=0$, i.e., on $\{p_h>0\}$, $d_h-1<0$. This observation achieves the proof of the lemma. \ok

\brem
The above lemma shows that, when the equation (\ref{equation-kappa}) has a solution, the first condition in Assumption \ref{1+fin} is necessary at least for $R\in \{S_m: 1\leq m<\mathsf{N}^i\}\cup \{T_n: 1\leq n<\mathsf{N}^a\}$. 
\erem

\bl\label{F-Sn}
For a fixed $1\leq n<\mathsf{N}^i$, let $(B_1,B_2,B_3,\ldots,B_{d})$ be a partition which satisfies the relation $
\mathcal{F}_{S_n}=\mathcal{F}_{S_n-}\vee\sigma(B_1,B_2,\ldots,B_{d}).
$
Denote $q_k=\mathbb{P}[B_k|\mathcal{F}_{S_n-}], 1\leq k\leq d$. Let $\mathtt{r}_n$ be the vector valued process introduced in Lemma \ref{varphi-r}.  
\ebe
\item[. ]
For any finite random variable $\xi\in\mathcal{F}_{S_n}$, the conditional expectation $\mathbb{E}[\xi|\mathcal{F}_{S_n-}]$ is well-defined. Let $$
\mathfrak{i}(\xi)_k=\ind_{\{q_k>0\}}\frac{1}{q_k}\mathbb{E}[\ind_{B_k}\xi|\mathcal{F}_{S_n-}], \ 1\leq k\leq d.
$$
We have $\xi=\sum_{h=1}^d\mathfrak{i}(\xi)_h\ind_{B_h}$. We can consider $\xi$ as a random function define on the measurable space $\{1,2,\ldots,d\}$ (the value of $\xi$ at point $k$ being $\mathfrak{i}(\xi)_k$).  
\item[. ]
Denote $n$-dimensional vector valued random variable $n_k=\mathfrak{i}(\Delta_{S_n}N)_k, 1\leq k\leq d$. We have $
(1+\transp\overline{\varphi}_{S_n}n_k)q_k
=
(1+\transp\overline{\varphi}_{S_n}\mathtt{R}_{n})\mathbb{E}[\ind_{B_k}|\mathcal{G}_{S_n-}]\geq 0
$
for $1\leq k\leq d$, on $\{S_n\leq T, S_n<\infty\}$.
  
\item[. ]
We have$$
(1+\transp\overline{\varphi}_{S_n}\mathtt{R}_{n})\mathbb{E}[(\!(v|\Delta_{S_n}W)\!)^2|\mathcal{G}_{S_n-}]
=
\sum_{h=0}^d(1+(\overline{\varphi}_{S_n}n_h))(\!(v|w_h)\!)^2q_h
$$
for $v\in\mathbb{R}^d$, $w_k=\mathfrak{i}(\Delta_{S_n}W)_k, 1\leq k\leq d$, on $\{S_n\leq T, S_n<\infty\}$.

\item[. ]
We have $(1+\transp\overline{\varphi}_{S_n}\mathtt{R}_{n})>0$ almost surely on $\{S_n\leq T, S_n<\infty\}$, and $(1+\transp\overline{\varphi}_{S_n}n_k)q_k>0$ if and only if $q_k>0$.

\item[. ] 
If $\chi_{S_n}$ denotes the random variable defined by $\chi_{S_n}=\sum_{k=1}^dk\ind_{B_k}\in\{1,2,\ldots,d\}$,
the conditional law of $\chi_{S_n}$ on $\{S_n\leq T, S_n<\infty\}$ given $\mathcal{G}_{S_n-}$ is absolutely continuous with respect to the conditional law of $\chi_{S_n}$ given $\mathcal{F}_{S_n-}$ with a (random) density function $\sum_{k=1}^d\frac{1+\transp\overline{\varphi}_{T_n}n_k}{1+\transp\overline{\varphi}_{S_n}\mathtt{R}_{n}}\ind_{\{k\}}$ 

\dbe

\el

\textbf{Proof.}
The proof of the first assertion is straightforward. To prove the second assertion, it is enough to notice that, in the same way as in Lemma \ref{conditional-expectation-i}, for any $1\leq k\leq d$, $$
\dcb
((1+\transp\overline{\varphi}_{S_n}\mathtt{R}_{n})\ind_{B_k}\ind_{[S_n,\infty)})^{\mathbb{G}-p}

&=&((1+\transp\overline{\varphi}_{S_n}n_k)q_k\ind_{[S_n,\infty)}
)^{\mathbb{G}-p}\\

\dce
$$
on $[0,T]$. To prove the third assertion, we note firstly the equality $$
(\!(v|\overline{\mathsf{ci}}_{S_n}g_{S_n} v)\!)
=
\transp\overline{\varphi}_{S_n}\mathbb{E}[\Delta_{S_n}N(\!(v|\Delta_{S_n}W)\!)^2|\mathcal{F}_{S_n-}], \ v\in\mathbb{R}^d,
$$
consequence of Lemma \ref{cig}. Applying then Lemma \ref{conditional-expectation-i}, we write, on $\{S_n\leq T, S_n<\infty\}$, $$
\dcb
&&(1+\transp\overline{\varphi}_{S_n}\mathtt{R}_{n})\mathbb{E}[(\!(v|\Delta_{S_n}W)\!)^2|\mathcal{G}_{S_n-}]\\

&=&\mathbb{E}[(\!(v|\Delta_{S_n}W)\!)^2|\mathcal{F}_{S_n-}]
+\transp\overline{\varphi}_{S_n}\mathbb{E}[\Delta_{S_n}N(\!(v|\Delta_{S_n}W)\!)^2|\mathcal{F}_{S_n-}]
\\

&=&\sum_{h=0}^d(1+\transp\overline{\varphi}_{S_n}n_h)(\!(v|w_h)\!)^2p_h
\dce
$$
Consider $(1+\transp\overline{\varphi}_{S_n}\mathtt{R}_{n})$. We compute on $[0,T]$ $$
\dcb
(\ind_{\{1+\transp\overline{\varphi}_{S_n}\mathtt{R}_{n}=0\}}\ind_{[S_n,\infty)})^{\mathbb{G}-p}

&=&\ind_{\{1+\transp\overline{\varphi}\mathtt{r}_{n}=0\}}(1+\transp\overline{\varphi}\mathtt{r}_{n})\centerdot(\ind_{[S_n,\infty)})^{\mathbb{F}-p}=0
\dce
$$
This yields $
\mathbb{E}[\ind_{\{1+\transp\overline{\varphi}_{S_n}\mathtt{R}_{n}=0\}}\ind_{\{S_n\leq T,S_n<\infty\}}]
=0,
$
proving the fourth assertion. The last assertion follows from the preceding ones. \ok

\bcor\label{varphi-N}
Let $R$ be a $\mathbb{F}$ stopping time. Suppose that $R$ is either predictable or totally inaccessible. Then, $
1+\transp\overline{\varphi}_R\Delta_R N> 0.
$
\ecor

\textbf{Proof.} We consider only the case where $R$ is predictable. We can suppose without loss of generality that $R$ is one of the $T_n, 1\leq n<\mathsf{N}^a$. We compute for $0\leq k\leq d$ :$$
\dcb
0\leq\mathbb{E}[\ind_{\{1+\transp\overline{\varphi}_Rn_k \leq 0\}}\ind_{A_k}|\mathcal{G}_{T_n-}]

&=&\ind_{\{1+\transp\overline{\varphi}_Rn_k \leq 0\}}(1+\transp\overline{\varphi}_{R}n_k)p_k\leq 0
\dce
$$
It follows that $\ind_{\{1+\transp\overline{\varphi}_R\Delta_R N \leq 0\}}\ind_{A_k}=0$ for $0\leq k\leq d$, i.e. $1+\transp\overline{\varphi}_R\Delta_R N>0$. \ok

\

\subsection{Consequences of Assumption \ref{1+fin}}

\bl\label{78}
Under Assumption \ref{1+fin}, the inequalities (\ref{di-u-inequality}) and (\ref{da-u-inequality}) are satisfied.
\el

\textbf{Proof.}
Let $\mathsf{u}$ be the process introduced in Assumption \ref{1+fin}. Let $1\leq n<\mathsf{N}^a,  1\leq m<\mathsf{N}^i$. We use the notation $n_\ast$ introduced in Lemma \ref{F-Tn} and Lemma \ref{F-Sn} representing the values of $\Delta N$ at $T_n$ or at $S_m$. Assumption \ref{1+fin} implies $$
\dcb
(1+\transp\overline{\varphi}_{T_n}n_k )\ind_{A_k}\geq \mathsf{u}_{T_n}\ind_{A_k},\  0\leq k\leq d, \mbox{ on $\{T_n\leq T, T_n<\infty\}$}\\
(1+\transp\overline{\varphi}_{S_m}n_j )\ind_{B_j}\geq \mathsf{u}_{S_m}\ind_{B_j},\ 1\leq j\leq d, \mbox{ on $\{S_m\leq T, S_m<\infty\}$}.
\dce
$$
Taking the conditional expectations $\mathbb{E}[\cdot|\mathcal{G}_{T_n-}]$ or $\mathbb{E}[\cdot|\mathcal{G}_{S_n-}]$ on the sets $\{T_n\leq T, T_n<\infty\}$ or respectively $\{S_m\leq T, S_m<\infty\}$, we obtain
$$
\dcb
(1+\transp\overline{\varphi}_{T_n}n_k )(1+\transp\overline{\varphi}_{T_n}n_k)p_k\geq \mathsf{u}_{T_n}(1+\transp\overline{\varphi}_{T_n}n_k)p_k,\\
(1+\transp\overline{\varphi}_{S_m}n_j )\frac{(1+\transp\overline{\varphi}_{S_m}n_j)}{(1+\transp\overline{\varphi}_{S_m}\mathtt{R}_{m})}q_j
\geq \mathsf{u}_{S_m}\frac{(1+\transp\overline{\varphi}_{S_m}n_j)}{(1+\transp\overline{\varphi}_{S_m}\mathtt{R}_{m})}q_j.
\dce
$$
Assumption \ref{1+fin} ensures $\{(1+\transp\overline{\varphi}_{T_n}n_k)p_k=0\}=\{p_k=0\}$, whilst Lemma \ref{F-Sn} implies $\{(1+\transp\overline{\varphi}_{S_m}n_j)q_j=0\}=\{q_j=0\}$. 
The above inequalities become$$
\dcb
(1+\transp\overline{\varphi}_{T_n}n_k )p_k\geq \mathsf{u}_{T_n}p_k\\
(1+\transp\overline{\varphi}_{S_m}n_j )q_j\geq \mathsf{u}_{S_m}q_j
\dce
$$ 
Under the same condition, for a $v\in\mathbb{R}^d$, according to Lemma \ref{at-Tn} and \ref{F-Tn}, $$
\dcb
&&(\!(v|(\Delta_{T_n}a)\transp(\mathfrak{I}_d+\overline{\mathsf{ce}}-\Delta a\overline{\mathsf{f}}\transp\overline{\mathsf{f}})_{T_n} v)\!)\\

&=&\sum_{h=0}^d(1+(\transp\overline{\varphi}_{T_n}n_h))p_h\left((\!(v|w_h)\!)-\sum_{k=0}^d(1+(\transp\overline{\varphi}_{T_n}n_k))(\!(v|w_k)\!)p_k\right)^2\\

&\geq&\mathsf{u}_{T_n}(\!(v|(\Delta_{T_n}a)\!) v)
\dce
$$
proving the inequality (\ref{da-u-inequality}), and by Lemma \ref{conditional-expectation-i} and \ref{F-Sn}
$$
\dcb
(\!(v|g_{S_n}\transp(\mathfrak{I}_d+\overline{\mathsf{ci}})_{S_n} v)\!)

=
\sum_{h=1}^d(1+\transp\overline{\varphi}_{S_m}n_h)(\!(v|w_h)\!)^2q_h

&\geq&\mathsf{u}_{S_m} (\!(v|g_{S_n} v)\!)
\dce
$$
proving the inequality (\ref{di-u-inequality}). \ok

\bl\label{integrable-xi}
Suppose Assumption \ref{1+fin}. Let $$
\dcb
\xi_n=\ind_{\{S_n\leq T, S_n<\infty\}}\ind_{\{g_{S_n}\neq 0\}}\left(\transp(\mathfrak{I}_d+\overline{\mathsf{ci}})_{S_n}\right)^*\mathtt{p}_{g_{S_n}}(\overline{d}+\overline{\mathsf{f}})_{S_n}\\
\xi'_n=\ind_{\{T_n\leq T, T_n<\infty\}}\ind_{\{\Delta_{T_n}a\neq 0\}}\left(\transp(\mathfrak{I}_d+\overline{\mathsf{ce}}-\Delta a\overline{\mathsf{f}}\transp\overline{\mathsf{f}})_{T_n}\right)^*\mathtt{p}_{\Delta_{T_n}a}(\overline{d}+\overline{\mathsf{f}})_{T_n}
\dce
$$ 
as in Lemma \ref{di-xi-lemma}, or respectively in Lemma \ref{da-xi-lemma}. Then, for any finite $\mathbb{G}$ stopping time $R\leq T$, $$
\dcb
\mathbb{E}[\sum_{1\leq n<\mathsf{N}^i}(\transp\xi_{n}\Delta_{S_n}W)^2\ind_{\{S_n\leq R\}}]
\leq
\mathbb{E}[\frac{1}{\mathsf{u}}\transp(\overline{d}+\transp\overline{\mathsf{f}})(\centerdot b)(\overline{d}+\transp\overline{\mathsf{f}})_R]\\

\mathbb{E}[\sum_{1\leq n<\mathsf{N}^a}(\transp\xi'_{n}\Delta_{T_n}\widetilde{W})^2\ind_{\{T_n\leq R\}}]
\leq
\mathbb{E}[\frac{1}{\mathsf{u}}\transp(\overline{d}+\transp\overline{\mathsf{f}})(\centerdot a)(\overline{d}+\transp\overline{\mathsf{f}})_R]
\dce
$$
\el

\textbf{Proof.} Note firstly that, by Lemma \ref{78}, $\xi_n$ and $\xi'_n$ are well-defined. By Lemma \ref{W-jumps}, $$
\dcb
\transp(\overline{d}+\transp\overline{\mathsf{f}})(\centerdot b)(\overline{d}+\transp\overline{\mathsf{f}})
\geq
\sum_{1\leq n<\mathsf{N}^i}
(\overline{d}+\transp\overline{\mathsf{f}})(\centerdot\left(g_{S_n}\ind_{[S_n,\infty)}\right)^{\mathbb{F}-p})(\overline{d}+\transp\overline{\mathsf{f}})\\

\transp(\overline{d}+\transp\overline{\mathsf{f}})(\centerdot a)(\overline{d}+\transp\overline{\mathsf{f}})
\geq
\sum_{1\leq n<\mathsf{N}^a}(\overline{d}+\transp\overline{\mathsf{f}})_{T_n}(\Delta_{T_n}a)(\overline{d}+\transp\overline{\mathsf{f}})_{T_n}\ind_{[T_n,\infty)}
\dce
$$
With this in mind, we see that the lemma is proved, once we prove the individual inequalities at $S_n$ or at $T_{n'}$.

According to Lemma \ref{78}, Assumption \ref{1+fin} implies the inequalities (\ref{di-u-inequality}) and (\ref{da-u-inequality}) with $\varkappa_n=\mathsf{u}_{T_n}$ and $\varrho_n=\mathsf{u}_{S_n}$. In Lemma \ref{78}, we proved $(1+\transp\overline{\varphi}_{S_n}n_j)q_j\geq \mathsf{u}_{S_n}q_j, 1\leq j\leq d$ on $\{S_m\leq T, S_m<\infty\}$. Summing up these inequalities over $1\leq j\leq d$, we prove $1+\transp\overline{\varphi}\mathtt{R}_{n}\geq \mathsf{u}_{S_n}$.
Applying Lemma \ref{conditional-expectation-i} and Lemma \ref{symmetric-matrix}, we obtain, on $\{S_m\leq T, S_m<\infty\}$:$$
\dcb
(1+\transp\overline{\varphi}_{S_n}\mathtt{R}_{n})\mathbb{E}[(\transp \xi_n\Delta_{S_n}W)^2|\mathcal{G}_{S_n-}]

&\leq&\frac{1}{\mathsf{u}_{S_n}}
(\!((\overline{d}+\overline{\mathsf{f}})_{S_n}|g_{S_n}(\overline{d}+\overline{\mathsf{f}})_{S_n})\!)\\
\dce
$$
Therefore, by Lemma \ref{F-Sn},
$$
\dcb
&&\mathbb{E}[(\transp \xi_n\Delta_{S_n}W)^2\ind_{\{S_n\leq R\}}]\\

&\leq&\mathbb{E}[\frac{1}{(1+\transp\overline{\varphi}_{S_n}\mathtt{R}_{n})}\ \frac{1}{\mathsf{u}_{S_n}}
(\!((\overline{d}+\overline{\mathsf{f}})_{S_n}|g_{S_n}(\overline{d}+\overline{\mathsf{f}})_{S_n})\!)\ind_{\{S_n\leq R\}}]\\

&=&\mathbb{E}[\frac{1}{\mathsf{u}}(\transp\overline{d}+\transp\overline{\mathsf{f}}) 
(\centerdot(g_{S_n}\ind_{[S_n,\infty)})^{\mathbb{F}-p})(\overline{d}+\overline{\mathsf{f}})_R]\\

\dce
$$
This proves the first inequality. The second inequality can be proved in the same way. \ok 

\

\subsection{Proof of the main theorem}

Suppose Assumptions \ref{assump0}, \ref{assump1} and \ref{assump-mrt}. Suppose Assumption \ref{1+fin} with the process $\mathsf{u}$. Suppose that the processes $
\frac{1}{\mathsf{u}}\centerdot[D^d,D^d] \ \mbox{ and \ } \frac{1}{\mathsf{u}}\transp\overline{\varphi}(\centerdot[N^d,\transp N^d])\overline{\varphi}
$
are $(\mathbb{P},\mathbb{G})$ locally integrable. According to Lemma \ref{78}, the inequalities (\ref{di-u-inequality}) and (\ref{da-u-inequality}) are satisfied. According to Lemma \ref{integrable-xi}, we have the following estimations on the random variables $\xi_n$ and $\xi'_n$ as defined in Lemma \ref{di-xi-lemma} and respectively in Lemma \ref{da-xi-lemma} :$$
\dcb
\mathbb{E}[\sum_{1\leq n<\mathsf{N}^i}(\transp\xi_{n}\Delta_{S_n}W)^2\ind_{\{S_n\leq R\}}]
\leq
\mathbb{E}[\frac{1}{\mathsf{u}}\transp(\overline{d}+\transp\overline{\mathsf{f}})(\centerdot b)(\overline{d}+\transp\overline{\mathsf{f}})_R]\\ 

\mathbb{E}[\sum_{1\leq n<\mathsf{N}^a}(\transp\xi'_{n}\Delta_{T_n}\widetilde{W})^2\ind_{\{T_n\leq R\}}]
\leq
\mathbb{E}[\frac{1}{\mathsf{u}}\transp(\overline{d}+\transp\overline{\mathsf{f}})(\centerdot a)(\overline{d}+\transp\overline{\mathsf{f}})_R]
\dce
$$
where $R$ is any finite $\mathbb{G}$ stopping time inferior or equal to $T$. We note that$$
\dcb
\sqrt{\mathbb{E}[\frac{1}{\mathsf{u}}\transp(\overline{d}+\transp\overline{\mathsf{f}})(\centerdot b)(\overline{d}+\transp\overline{\mathsf{f}})_R]}
&\leq &
\sqrt{\mathbb{E}[\frac{1}{\mathsf{u}}\transp\overline{d}(\centerdot b)\overline{d}_R]}
+
\sqrt{\mathbb{E}[\frac{1}{\mathsf{u}}\transp\overline{\mathsf{f}}(\centerdot b)\transp\overline{\mathsf{f}}_R]}\\

&= &
\sqrt{\mathbb{E}[\frac{1}{\mathsf{u}}\centerdot [D^{di},D^{di}]^{\mathbb{F}-p}_R]}
+
\sqrt{\mathbb{E}[\frac{1}{\mathsf{u}}\transp\overline{\varphi}(\centerdot[N^{di},\transp N^{di}])\overline{\varphi}_R]}\\

&\leq&
\sqrt{\mathbb{E}[\frac{1}{\mathsf{u}}\centerdot [D^d,D^d]^{\mathbb{F}-p}_R]}
+
\sqrt{\mathbb{E}[\frac{1}{\mathsf{u}}\transp\overline{\varphi}(\centerdot[N^d,\transp N^d])\overline{\varphi}_R]}\\

\dce
$$ 
In the same way we obtain $$
\sqrt{\mathbb{E}[\frac{1}{\mathsf{u}}\transp(\overline{d}+\transp\overline{\mathsf{f}})(\centerdot a)(\overline{d}+\transp\overline{\mathsf{f}})_R]}
\leq
\sqrt{\mathbb{E}[\frac{1}{\mathsf{u}}\centerdot [D^d,D^d]^{\mathbb{F}-p}_R]}
+
\sqrt{\mathbb{E}[\frac{1}{\mathsf{u}}\transp\overline{\varphi}(\centerdot[N^d,\transp N^d])\overline{\varphi}_R]}
$$
It follows that the increasing processes $$
\sum_{1\leq n<\mathsf{N}^i}(\transp\xi_{n}\Delta_{S_n}W)^2\ind_{[S_n,\infty)},\ \
\sum_{1\leq n<\mathsf{N}^a}(\transp\xi'_{n}\Delta_{T_n}\widetilde{W})^2\ind_{[T_n,\infty)}
$$
are $(\mathbb{P},\mathbb{G})$ locally integrable. Lemma \ref{da-xi-lemma} and \ref{di-xi-lemma} are valid. Let $\overline{K}$ and $\overline{K}'$ be the $\mathbb{G}$ predictable processes defined in Lemma \ref{di-xi-lemma} and respectively in Lemma \ref{da-xi-lemma}. Then, $\overline{K}$ solves the equation (\ref{equation-rho}), whilst $\overline{K}'$ solves the equation (\ref{equation-kappa}). Now in order to solve the structure condition (\ref{structure-condition-di}) and (\ref{structure-condition-da}), we need only to show that $\transp\overline{K}\Delta\widetilde{W^{di}}<1$ and $\transp\overline{K}'\Delta\widetilde{W^{da}}<1$.

We will consider only the first inequality. The second one can be dealt with similarly. By Lemma \ref{di-xi-lemma}, $$
\transp\overline{K}\centerdot \widetilde{W^{di}}
=\sum_{1\leq n<\mathsf{N}^i}(\xi_{n} \Delta_{S_n}W\ind_{[S_n,\infty)}-(\xi_{n}\Delta_{S_n}W\ind_{[S_n,\infty)})^{\mathbb{G}-p})
$$
on $[0,T]$. By Lemma \ref{varphi-r}, $(\xi_{n}\Delta_{S_n}W\ind_{[S_n,\infty)})^{\mathbb{G}-p})$ is continuous. To show $\transp\overline{K}\Delta\widetilde{W^{di}}<1$, it is enough to show $\xi_{n}\Delta_{S_n}W<1$ for any $1\leq n<\mathsf{N}^i$. In the proof of Lemma \ref{di-xi-lemma}, we obtained  
$$
(1+\transp\overline{\varphi}_{S_n}\mathtt{R}_n)\transp\xi_n\mathbb{E}[\Delta_{S_n}W\transp \Delta_{S_n}W|\mathcal{G}_{S_n-}]
=(\transp\overline{d}+\transp\overline{\mathsf{f}})_{S_n}\mathbb{E}[\Delta_{S_n}W\transp \Delta_{S_n}W|\mathcal{F}_{S_n-}]
$$
on $\{S_n\leq T, S_n<\infty\}$. Multiplying this identity by $\ddag_{S_n}\ind_{B_k}$ for a $1\leq k\leq d$, we get
$$
(1+\transp\overline{\varphi}_{S_n}\mathtt{R}_n)\transp\xi_n\mathbb{E}[\Delta_{S_n}W \ind_{B_k}|\mathcal{G}_{S_n-}]
=(\transp\overline{d}+\transp\overline{\mathsf{f}})_{S_n}\mathbb{E}[\Delta_{S_n}W \ind_{B_k}|\mathcal{F}_{S_n-}]
$$
This yields, with the notations in Lemma \ref{F-Sn},$$
\transp\xi_nw_k(1+\transp\overline{\varphi}_{S_n}n_k)q_k
=(\transp\overline{d}+\transp\overline{\mathsf{f}})_{S_n}w_kq_k
$$
Either, $w_kq_k=0$ in which case $\Delta_{S_n}W\ind_{B_k}=0$ and consequently $\xi_{n}\Delta_{S_n}W\ind_{B_k}<1$. Either $w_kq_k\neq 0$ in which case the above identity implies$$
\transp\xi_nw_k(1+\transp\overline{\varphi}_{S_n}n_k)\ind_{B_k}
=(\transp\overline{d}+\transp\overline{\mathsf{f}})_{S_n}w_k\ind_{B_k}
$$
or equivalently
$$
\transp\xi_n\Delta_{S_n}W(1+\transp\overline{\varphi}_{S_n}\Delta_{S_n}N)\ind_{B_k}
=(\Delta_{S_n}D+\transp\overline{\varphi}_{S_n}\Delta_{S_n}N)\ind_{B_k}
$$
By Corollary \ref{varphi-N}, $(1+\transp\overline{\varphi}_{S_n}\Delta_{S_n}N)>0$ on $B_k$. The above identity becomes
$$
\transp\xi_n\Delta_{S_n}W\ind_{B_k}
=\frac{\Delta_{S_n}D+\transp\overline{\varphi}_{S_n}\Delta_{S_n}N}{1+\transp\overline{\varphi}_{S_n}\Delta_{S_n}N}\ind_{B_k}
$$
By Assumption \ref{assump0}, $\Delta D< 1$. We conclude finally $\transp\xi_n\Delta_{S_n}W\ind_{B_k}<1$. This being valid for any $1\leq n<\mathsf{N}^i$ and $1\leq k\leq d$, we prove the first inequality.

Theorem \ref{main} is proved.

\

\end{document}